\documentclass{article}

\usepackage{amsmath}
\usepackage{latexsym}
\usepackage{graphicx}
\usepackage{psfrag}
\usepackage{amssymb}
\usepackage{amscd}
\usepackage{fullpage}

\newcommand{\comment}[1]{}

  \newtheorem{thm}{Theorem}[section]
  \newtheorem{defn}[thm]{Definition}
  \newtheorem{lemma}[thm]{Lemma}
  \newtheorem{cor}[thm]{Corollary}
  \newtheorem{prop}[thm]{Proposition}

  \newtheorem{remark}[thm]{Remark}

\def\proof{{\noindent{\it Proof.\ }}}
\def\endproof{{\hfill $\Box$}}

\def\<{\langle}
\def\>{\rangle}
\def\0{{{\bf 0}}}

\def\CA{{\mathcal A}}
\def\CM{{\mathcal M}}
\def\CX{{\mathcal X}}

\def\bq{{\underline q}}

\def\CC{{\mathbb C}}

\def\FF{{\mathbb F}}
\def\GG{{\mathbb G}}

\def\QQ{{\mathbb Q}}

\def\TT{{\mathbb T}}
\def\ZZ{{\mathbb Z}}

\def\OO{{\mathcal O}}

\def\mm{{\mathfrak m}}
\def\pp{{\mathfrak p}}
\def\ll{{\mathfrak l}}
\def\tJ{{\tilde J}}
\def\tx{{\tilde x}}
\def\tmm{{\tilde {\mathfrak m}}}
\def\tTT{{\tilde {\mathbb T}}}
\def\Dp{{D^{\prime}}}
\def\tA{{\tilde {\mathcal A}}}
\def\dual{\dag}

\newcommand{\spec}{\operatorname{Spec}}

\newcommand{\tr}{\operatorname{Tr}}
\newcommand{\Hom}{\operatorname{Hom}}
\newcommand{\Id}{\operatorname{Id}}

\newcommand{\End}{\operatorname{End}}
\newcommand{\Frob}{\operatorname{Frob}}
\newcommand{\Ver}{\operatorname{Ver}}
\newcommand{\gal}{\operatorname{Gal}}
\newcommand{\Aut}{\operatorname{Aut}}
\newcommand{\GL}{\operatorname{GL}}
\newcommand{\SL}{\operatorname{SL}}

\newcommand{\ann}{\operatorname{Ann}}
\newcommand{\Rnew}{\mbox{\rm \tiny $R$-new}}
\newcommand{\Unew}{\mbox{\rm \tiny $U$-new}}
\newcommand{\Uold}{\mbox{\rm \tiny $U$-old}}

\newcommand{\pnew}{\mbox{\rm \tiny $p$-new}}
\newcommand{\pold}{\mbox{\rm \tiny $p$-old}}
\newcommand{\old}{\mbox{\rm \tiny old}}
\newcommand{\new}{\mbox{\rm \tiny new}}
\newcommand{\qold}{\mbox{\rm \tiny $q$-old}}
\newcommand{\qnew}{\mbox{\rm \tiny $q$-new}}
\newcommand{\pqnew}{\mbox{\rm \tiny $pq$-new}}
\newcommand{\ppnew}{\mbox{\rm \tiny $pp^{\prime}$-new}}
\newcommand{\qqnew}{\mbox{\rm \tiny $q_1q_2$-new}}
\newcommand{\qtold}{\mbox{\rm \tiny $q_2$-old}}

\newcommand{\NDnew}{\mbox{\rm \tiny $\frac{N}{D}$-new}}
\newcommand{\RDnew}{\mbox{\rm \tiny $\frac{R}{D}$-new}}
\newcommand{\Jmin}{J^{\mbox{\rm \tiny min}}}
\newcommand{\tJmin}{{\tilde J}^{\mbox{\rm \tiny min}}} 

\newcommand{\red}{\mbox{\rm \tiny red}}

	{\begin{list}
		{\noindent\makebox[0mm][r]{\arabic{enumi}.}}
		{\usecounter{enumi} \topsep=1.5mm \itemsep=0mm}
	}
	{\end{list}}

	{\begin{list}
		{\noindent\makebox[0mm][r]{\arabic{enumi}.}}
		{\usecounter{enumi} \topsep=1.5mm \itemsep=-1mm}
	}
	{\end{list}}

\begin{document}

\title{On maps between modular Jacobians  
and Jacobians of Shimura curves}
\author{David Helm\footnote{
Harvard University; dhelm@math.harvard.edu}}

\maketitle

Fix a squarefree integer $N$, divisible by an even number of primes,
and let $\Gamma^{\prime}$ be a congruence subgroup of level $M$,
where $M$ is prime to $N$.   For each $D$ dividing $N$ and divisible
by an even number of primes, 
the Shimura curve $X^D(\Gamma_0(N/D) \cap \Gamma^{\prime})$ associated
to the indefinite quaternion algebra of discriminant $D$ and 
$\Gamma_0(N/D) \cap \Gamma^{\prime}$-level
structure is well defined, and we can consider its Jacobian
$J^D(\Gamma_0(N/D) \cap \Gamma^{\prime})$.  Let $J^D$ denote the
$\frac{N}{D}$-new subvariety of this Jacobian. 

By the Jacquet-Langlands correspondence~\cite{JL} and Faltings' 
isogeny theorem~\cite{Fa}, there are Hecke-equivariant among 
the various varieties $J^D$ defined above.  However, since the isomorphism of 
Jacquet-Langlands is noncanonical, this perspective gives no
information about the isogenies so obtained beyond their 
existence.  In this paper, we study maps between the varieties $J^D$ in
terms of the maps they induce on the character groups of the tori
corresponding to the mod $p$ reductions of these varieties for $p$
dividing $N$. 
Our characterization of such maps in these terms allows us to classify
the possible kernels of maps from $J^D$ to $J^{D^{\prime}}$, 
for $D$ dividing $D^{\prime}$, up to support on a small finite
set of maximal ideals of the Hecke algebra.  This allows us to compute
the Tate modules $T_{\mm} J^D$ of $J^D$ at all non-Eisenstein $\mm$
of residue characteristic $l > 3.$  
These computations have implications for the multiplicities
of irreducible Galois representations in the torsion of Jacobians of Shimura
curves; one such consequence 
is a ``multiplicity one'' result for Jacobians of Shimura
curves.

\section{Introduction}

The Jacquet-Langlands correspondence asserts the existence of noncanonical
isomorphisms between spaces of modular forms coming from different
quaternion algebras over $\QQ$.  More precisely, if $D$ is a squarefree
product of an even number of primes, $p$ and $q$ are distinct primes
not dividing $D$, and $B$ and $B^{\prime}$ are quaternion algebras over
$\QQ$ with discriminants $D$ and $Dpq$, respectively, then the 
Jacquet-Langlands correspondence asserts the existence of an isomorphism:
$$S_2^{Dpq}(\Gamma) \cong S_2^D(\Gamma_0(pq) \cap \Gamma)_{\pqnew}.$$
Here $S_2^{Dpq}(\Gamma)$ denotes the space of weight two modular forms 
associated to the quaternion algebra $B^{\prime}$ and the congruence
subgroup $\Gamma$.  Similarly, $S_2^D(\Gamma_0(pq) \cap \Gamma)$ is
the space of weight two modular forms associated to the quaternion
algebra $B$ and the congruence subgroup $\Gamma_0(pq) \cap \Gamma$.
The subscript $pq$-new denotes the subspace of 
$S_2^D(\Gamma_0(pq) \cap \Gamma)$ consisting of forms which are ``new''
for both $p$ and $q$. 

The existence of such an isomorphism has consequences for the geometry
of modular curves and Shimura curves.  In particular, let
$X^{Dpq}(\Gamma)$ denote the Shimura curve associated to
the quaternion algebra $B^{\prime}$ and the congruence subgroup 
$\Gamma_0(pq) \cap \Gamma$, and define $X^D(\Gamma_0(pq) \cap \Gamma)$
similarly.
Then we may consider their Jacobians $J^D(\Gamma_0(pq) \cap \Gamma)$
and $J^{Dpq}(\Gamma)$.  The Jacquet-Langlands correspondence,
together with Eichler-Shimura theory, tells us that, after tensoring
with $\QQ$, the $l$-adic Tate module
of $J^{Dpq}(\Gamma)$ is isomorphic to the corresponding
Tate module of the $pq$-new subvariety
of $J^D(\Gamma_0(pq) \cap \Gamma).$  It follows from~\cite{ribisog}
or from Faltings' isogeny theorem~\cite{Fa}, that $J^{Dpq}$ is isogenous
to the $pq$-new subvariety of $J^D(\Gamma_0(pq) \cap \Gamma)$.

Actually producing an isogeny between these two abelian varieties, on the other
hand, seems much more difficult.  The above argument, while enough
to establish the existence of such an isogeny, provides no hint as to
how an explicit map might be constructed.  

The first results to give geometric information about
the relationship between $J^{Dpq}(\Gamma)$ and 
$J^D(\Gamma_0(pq) \cap \Gamma)$ 
are due to Ribet~\cite{level}, who established an isomorphism 
between the character group of the torus associated to the mod $p$ 
reduction of $J^{Dpq}(\Gamma)$ and a certain subgroup of the character
group associated to the mod $q$ reduction of $J^D(\Gamma_0(pq) \cap \Gamma)$.  
This isomorphism played a key role in his proof of the level-lowering
theorem.  It also allowed him to show~\cite{Israel} 
that even in the case where $D=1$ and $\Gamma$ is trivial,
$J_0(pq)_{\pqnew}$ and $J^{pq}(1)$ could have decidedly different
structures.  In particular, he constructed examples of maximal ideals
$\mm$ of the Hecke algebra $\TT$ 
for which $J_0(pq)_{\pqnew}[\mm]$ had dimension two
yet $J^{pq}(1)[\mm]$ had dimension four.  In particular, this would happen
whenever the (mod $l$) Galois representation $\rho_{\mm}$ associated 
to $\mm$ was unramified at $q$ and satisfied $\rho_{\mm}(\Frob_q) = \pm 1.$
Thus, any isogeny between the two varieties would encode data about
the restriction of $\rho_{\mm}$ to decomposition groups at $p$ and $q$,
for every maximal ideal $\mm$ of $\TT$.  In particular, the structure 
of any isogeny between the two varieties was necessarily complicated.

In this paper we study isogenies between the varieties described above
using an approach first suggested by Mazur.  Instead of attempting to
produce isogenies between the two varieties directly, we study the entire
category of abelian varieties with an action of the Hecke-algebra which
are Hecke-equivariantly isogenous to $J^{Dpq}(\Gamma)$. 

More precisely, fix a
squarefree integer $N$, prime to the level $M$ of $\Gamma$, and for each $D$ 
dividing $N$ and divisible by an even number of primes, let
$J^D = J^D(\Gamma_0(N/D) \cap \Gamma)_{\NDnew}$. 
Let $\TT$ be the subalgebra of
the endomorphism ring of some such $J^D$
generated by the Hecke operators
$T_l$ for $l$ prime, and the diamond 
bracket operators $\<d\>$.  Repeated applications
of the Jacquet-Langlands correspondence show that $\TT$ is
independent (up to canonical isomorphism)
of the choice of $J^D$ used to define it, and
therefore acts simultaneously on all the $J^D$.  Moreover, the
abelian varieties $J^D$ are all 
$\TT$-equivariantly isogenous.  We let $\CA$ denote the category of abelian
varieties with a $\TT$-action which are $\TT$-equivariantly isogenous to
$J$, and develop a formalism which provides a concrete description of
$\CA$ in terms of $\TT$-modules. 

In particular, we show that if one is
willing to work ``up to support on $S$'', where $S$ is the set of
Eisenstein primes of $\TT$, that one can classify all abelian varieties
in $\CA$ in terms of isomorphism
classes of $\TT$-modules.  Specifically, to any abelian variety $J$
isogenous to $J^1$ one can associate a certain $\TT$-module $[J]$
which is free over $\ZZ$ and satisfies $[J] \otimes \QQ \cong
\TT \otimes \QQ$.  (Such modules were first introduced by Mazur~\cite{MazurEis},
who called them ``rank one'' modules.)
After formally inverting isogenies supported on $S$,
and maps of $\TT$-modules whose cokernels are supported on $S$,
this
functor becomes an antiequivalence of categories  (Theorem~\ref{thm:equiv}).
Sections~\ref{sec:invert} and~\ref{sec:isogeny} are devoted to making this
precise.
Moreover, if $J$ and $J^{\prime}$ are isogenous to $J^1$,
knowing the isomorphism classes of $[J]$ and $[J^{\prime}]$
as $\TT$-modules is often enough (Proposition~\ref{prop:JJ'}) to determine
the kernels of maps $J \rightarrow J^{\prime}$
up to support on $S$.  Thus computing $[J]$ and $[J^{\prime}]$ allows
us to construct the isomorphism class of $J^{\prime}$ in terms of the
isomorphism class of $J$.  In particular one can express the Tate modules
of $J^{\prime}$ at maximal ideals outside $S$ in terms of the Tate modules
of $J$ and the modules $[J]$ and $[J^{\prime}]$.

Moreover, we obtain a relationship between the character groups $X_p(J)$
of the tori coming from the mod $p$ reductions of $J$ for $p$ dividing $N$, 
and the module $[J]$, for any $J$ isogenous to $J^1$
(Proposition~\ref{prop:charisogeny}).
In Section~\ref{sec:chargp} we establish a key relationship between
the character groups of $J^D$ as $D$ varies.  These two results allow us
to compute the modules $[J^D]$ explicitly.  The key result we obtain is
Theorem~\ref{thm:main}, in which we obtain $[J^D]$ for $D=N$ in terms
of a certain module of maps between character groups.   

Sections~\ref{sec:maps} and~\ref{sec:raise} are devoted to the proof
of Theorem~\ref{thm:main}.  In section~\ref{sec:maps}, we introduce a
condition, called {\em controllability}, associated to a prime $p$
dividing $N$ and a maximal ideal $\mm$ of $\TT$.  Controllability of
$\mm$ at $p$ implies that the character group of a certain abelian
variety isogenous to $J^D$ becomes free of rank one after completing
at $\mm$ (c.f. Lemma~\ref{lemma:control}); 
in section~\ref{sec:maps} we exploit this fact to show that 
Theorem~\ref{thm:main} holds locally at any maximal ideal $\mm$ which is
controllable at some prime dividing $D$.

It remains to handle those maximal ideals which are not controllable at
any primes dividing $D$, which we do in Section~\ref{sec:raise}.
Although such ideals are not at all common, they are much more difficult
to handle, and the arguments presented here are of a very technical
sort.  (The reader who is primarily interested in our results, or who is 
willing to take this extreme case of Theorem~\ref{thm:main} on faith,
is encouraged to read Section~\ref{sec:main} first; once Theorem~\ref{thm:main}
is established we have no further need of the results of 
Section~\ref{sec:raise}.) 
Our basic strategy for dealing with these maximal ideals is a 
level-raising argument.  In particular, we fix such a maximal
ideal $\mm$, and introduce 
two primes $q_1$ and $q_2$ to the level, in such a way that $\mm$ is
controllable at both $q_1$ and $q_2$.  We establish an analogue
of Theorem~\ref{thm:main} in this new setting, using controllability
at $q_1$ and $q_2$, and then show that it implies the desired result. 

In Section~\ref{sec:main} we explore the consequences of
Theorem~\ref{thm:main}.  In particular, we apply Theorem~\ref{thm:main}
along with Proposition~\ref{prop:chargp} to obtain several expressions
for $[J^D]$ in terms of $[J^{D^{\prime}}]$ whenever $D^{\prime}$ divides $D$
(Theorem~\ref{thm:global1}, Theorem~\ref{thm:globalmult1}, and
Corollary~\ref{cor:homDD'}).  In particular, if $p$ is a divisor of
$D^{\prime}$, we find that the $\TT$-modules
$\Hom(J^{D^{\prime}}, J^D)$ and $\Hom(X_p(J^D), X_p(J^{D^{\prime}}))$
are isomorphic up to support on $S$, via the natural map which takes
a morphism of varieties to the map it induces on character groups. 
Via Proposition~\ref{prop:JJ'}, this allows us to determine
the kernels of all morphisms $J^D \rightarrow J^{D^{\prime}}$, up
to support on $S$.  We give several different descriptions of these kernels
(Corollaries~\ref{cor:global1},~\ref{cor:globalmult1}, and~\ref{cor:homtwo}) 
with varying degrees of explicitness.  The latter two expressions, for
instance, are in terms of character groups which can be computed
algorithmically~\cite{Kohel}.  The upshot of this is that given
the isomorphism class of $J^1$, one can explicitly compute the
isomorphism class of $J^D$ for any $D$, up to support on $S$. 

This recipe for constructing $J^D$ from $J^1$ allows us to compute
many invariants of $J^D$ in terms of $J^1$.  In particular one can
compute the $\mm$-adic Tate modules of $J^D$ in this fashion, for any
$\mm$ outside $S$.  In particular one obtains expressions for the
dimension of $J^D[\mm]$ in this manner.  In this fashion, we obtain a
bound on the dimension of $J^D[\mm]$ in terms of the number of primes
dividing $D$ at which $\mm$ is not controllable, a significant strengthening
of an earlier result of L. Yang~\cite{Yang}.  An important special case
of this bound is when $\mm$ is controllable at every prime and
$J^1[\mm]$ has dimension two; in this case $J^D[\mm]$ has dimension two
as well-- that is, we have ``multiplicity one'' for $J^D$ at $\mm$.

\section{Notation, background, and conventions}

Before we begin, we fix notation, explicitly state the background results
referred to in the introduction, and establish some conventions which will
be in use throughout the paper.

\begin{defn}  Let $N$ be an integer, and $H$ a subgroup
of $(\ZZ/N\ZZ)^{\times}$.  We define $\Gamma_H(N)$ to be the
subgroup of $\SL_2(\ZZ)$ consisting
of all matrices 
$\left( \begin{array}{cc} a & b \\ c & d \end{array} \right)$ 
in $\SL_2(\ZZ)$ for
which $c \equiv 0$ (mod $N$), $a$ and $d$ lie in a given subgroup $H$
of $(\ZZ/N\ZZ)^{\times}$ modulo $N$, and
$b$ is arbitrary.  By a congruence subgroup of level $N$, we mean
a subgroup of the form $\Gamma_H(N)$ for some $H$.
\end{defn}

Let $D$ be a squarefree integer divisible by an even number
of primes, $N$ an integer prime to $D$, and $\Gamma$ a congruence
subgroup of level $N$.  Then we denote by $X^D(\Gamma)$ the
Shimura curve associated to the indefinite quaternion algebra over
$\QQ$ of discriminant $D$, with $\Gamma$-level structure.  We let
$J^D(\Gamma)$ denote its Jacobian.  (If $D=1$ we take $X^D(\Gamma)$
to be the modular curve with $\Gamma$-level structure, instead.)
 
Let $T$ be a squarefree integer dividing $N$, and suppose
that $\Gamma$ has the form $\Gamma_0(T) \cap \Gamma^{\prime},$
where $\Gamma^{\prime}$ is a congruence subgroup of level $N/T$
and $\Gamma_0(T) \subset \SL_2(\ZZ)$ is the subgroup of matrices
$\left( \begin{array}{cc} a & b \\ c & d \end{array} \right)$ 
for which $c$ is divisible by $T$.
Let $p$ divide $T$.  Then there are two natural
degeneracy maps: 
$$\alpha_p, \beta_p: X^D(\Gamma) \rightarrow 
X^D(\Gamma_0(T/p) \cap \Gamma^{\prime})$$
These induce maps
$$\begin{array}{ccc}
J^D(\Gamma)^2 \rightarrow J^D(\Gamma_0(T/p) \cap \Gamma^{\prime}) &
\text{and} &
J^D(\Gamma_0(T/p) \cap \Gamma^{\prime}) \rightarrow J^D(\Gamma)
\end{array}$$
by Picard and Albanese functoriality, respectively.  We let
$J^D(\Gamma)_{\pold}$ denote the subvariety of $J^D(\Gamma)$ generated
by the images of the former maps, and 
$J^D(\Gamma)_{\pnew}$ denote the connected component of the intersection
of the kernels of the latter maps.
Both are abelian subvarieties of $J^D(\Gamma)$, stable under the
action of the Hecke algebra.  If $U$ is a divisor of $T$, we
let $J^D(\Gamma)_{\Uold}$ (resp. $J^D(\Gamma)_{\Unew}$) denote
the connected components of the intersections of the
varieties $J^D(\Gamma)_{\pold}$ (resp. $J^D(\Gamma)_{\pnew}$)
for $p$ dividing $U$. 

We denote by $\TT^D(\Gamma)$ the subalgebra of $\End(J^D(\Gamma))$
generated by the Hecke operators $T_l$ for $l$ prime to $N$,
and the diamond bracket operators
$\<d\>$ for $d$ in $(\ZZ/N\ZZ)^{\times}$.  (When $l$ divides $N$ we
take the $T_l$ to be defined as in~\cite{level}.)  The
algebra $\TT^D(\Gamma)$ is a finite, flat
$\ZZ$-algebra.  If $\Gamma = \Gamma_0(T) \cap \Gamma^{\prime}$ as
above, and $U$ divides $\TT$, then $\TT^D(\Gamma)$ acts on
$J^D(\Gamma)_{\Unew}$ (resp. $J^D(\Gamma)_{\Uold}$) through a finite
quotient which is also flat over $\ZZ$.  We denote this
quotient by $\TT^D(\Gamma)^{\Unew}$ (resp. $\TT^D(\Gamma)^{\Uold}$).

If $J$ is an abelian variety with a $\TT^D(\Gamma)$-action,
then this action allows us to speak of the $U$-new subvarieties and quotients
of $J$ as well.  We take the $U$-new subvariety of $J$ to be the
connected component of the subvariety of $J$ annihilated by the
kernel $I$ of the map $\TT^D(\Gamma) \rightarrow \TT^D(\Gamma)^{\Unew}$,
and the $U$-new quotient of $J$ to be $J/IJ$.  The $U$-old subvariety
and quotient of $J$ are constructed in a similar manner.  It is 
straightforward to check that this construction agrees with the previous
definition of $J^D(\Gamma)^{\Unew}$ and $J^D(\Gamma)_{\Unew}$. 

We have the following relationship between the subvarieties and quotients
defined above and duality of abelian varieties:
\begin{prop}
Suppose $J$ is isogenous to $J^D(\Gamma)$.  Then
have $(J_{\Unew})^{\vee} \cong (J^{\vee})^{\Unew}.$
\end{prop}
\proof Consider the exact sequence:
$$0 \rightarrow J_{\Unew} \rightarrow J \rightarrow A \rightarrow 0.$$
Taking duals, we obtain an exact sequence:
$$0 \rightarrow A^{\vee} \rightarrow J^{\vee} \rightarrow (J_{\Unew})^{\vee}
\rightarrow 0.$$
(Here the action of $\TT$ on $J^{\vee}$ is induced by duality of abelian
varieties.)
It thus suffices to show that we can identify $A^{\vee}$ with
$IJ^{\vee}$, where $I$ is the kernel of the map
$\TT^D(\Gamma) \rightarrow \TT^D(\Gamma)^{\Unew}$.  Since $I$ annihilates 
$(J_{\Unew})^{\vee}$, we must have $IJ^{\vee} \subset A$. 
Counting ranks of abelian varieties, we find that the two must be equal.
\endproof

Fix a particular $D$ and $N$, a congruence subgroup
$\Gamma^{\prime}$ of level $N$, and two distinct primes $p$ and $q$
which do not divide $DN$.  Let $\Gamma = \Gamma^{\prime} \cap \Gamma_0(pq)$.
The tangent spaces $S_2^D(\Gamma)_{\pqnew}$ and $S_2^{Dpq}(\Gamma^{\prime})$
of $J^D(\Gamma)_{\pqnew}$ and $J^{Dpq}(\Gamma^{\prime})$ over $\QQ$, 
respectively, are modules for $\TT^D(\Gamma)^{\pqnew}$ and 
$\TT^{Dpq}(\Gamma^{\prime})$.

\begin{thm} \label{thm:J-L} (Jacquet-Langlands, Ribet)
There is a natural isomorphism: $\TT^D(\Gamma)^{\pqnew} \rightarrow 
\TT^{Dpq}(\Gamma^{\prime})$
which sends an operator $T_l$
in one algebra to the corresponding operator $T_l$ in the other (and
similarly for $\<d\>$.) 
Moreover, the $\TT^{Dpq}(\Gamma^{\prime})$-modules
$S_2^D(\Gamma)_{\pqnew}$ and $S_2^{Dpq}(\Gamma^{\prime})$ are
isomorphic.
\end{thm}
\proof Jacquet and Langlands~\cite{JL} showed this for the subalgebras of
$\TT^D(\Gamma)^{\pqnew}$ and $\TT^{Dpq}(\Gamma^{\prime})$ generated
by the diamond bracket operators and the Hecke operators $T_l$ for
$l$ prime to $DpqN$.  The results of Ribet~\cite{level} imply that
it holds for the full Hecke algebra.
\endproof

\begin{cor}
There is a $\TT$-equivariant isogeny:
$$J^D(\Gamma)_{\pqnew} \rightarrow J^{Dpq}(\Gamma^{\prime}).$$
\end{cor}
\proof Let $\TT = \TT^{Dpq}(\Gamma^{\prime})$, and $V_l$ and $V_l^{\prime}$
the $l$-adic Tate modules of $J= J^D(\Gamma)_{\pqnew}$ and 
$J^{\prime} = J^{Dpq}(\Gamma^{\prime})$, respectively.  Then $V_l \otimes \QQ$ and
$V_l^{\prime} \otimes \QQ$ are both free of rank two over $\TT \otimes \QQ$,
and for any $l$ prime to $DpqN$, the trace and determinant of $\Frob_l$ on 
both modules is given by $T_l$ and $l\<l\>$, respectively.  It follows
that $V_l \otimes \QQ$ and $V_l^{\prime} \otimes \QQ$ are $\TT$-equivariantly
isomorphic as $G_{\QQ}$-modules.  

By~\cite{Fa}, Theorem 4, 
the natural map:
$$\Hom(J,J^{\prime}) \otimes \QQ_l \rightarrow \Hom(V_l^{\prime},V_l) \otimes
\QQ_l$$
is an isomorphism.  This natural map takes $\TT$-equivariant morphisms to
$\TT$-equivariant morphisms, and hence induces an isomorphism:
$$\Hom_{\TT}(J,J^{\prime}) \otimes \QQ_l \rightarrow
\Hom_{\TT}(V_l^{\prime}, V_l) \otimes \QQ_l.$$  Thus $J$ and $J^{\prime}$
are $\TT$-equivariantly isogenous. 
\endproof 

We finish this section with some notation and conventions.

We will often be considering modules over a finite, flat $\ZZ$-algebra
$\TT$.  If $M$ is a $\TT$-module, and $\mm$ is a maximal ideal of $\TT$, 
then
$M_{\mm}$ denotes the {\em completion} of $M$ at $\mm$.  If $M$ is
free as a $\ZZ$-module, then $M^*$ denotes the $\ZZ$-dual $\Hom_{\ZZ}(M,\ZZ)$,
with the $\TT$-module structure induced from $M$.  If $M$ is free as
a $\ZZ_l$-module for some $l$, on the other hand, then by abuse of
notation we write $M^*$ for $\Hom_{\ZZ_l}(M,\ZZ_l).$  The advantage of
this notation is that $(M_{\mm})^*$ is naturally isomorphic to 
$(M^*)_{\mm}$ for any maximal ideal $\mm$ of $\TT$. 

The $\TT$-modules we are most interested in are all torsion free over $\ZZ$.
Unfortunately, the tensor product of two $\TT$-modules without $\ZZ$-torsion
often has $\ZZ$-torsion, which we wish to ignore.  Therefore we adopt
the convention that if $M$ and $M^{\prime}$ are torsion free $\TT$-modules
or $\TT_{\mm}$-modules, then $M \otimes M^{\prime}$ denotes their
tensor product (over $\TT$ or $\TT_{\mm}$) {\em modulo $\ZZ$-torsion}.
It can be verified that this operation satisfies the universal property
of the tensor product in the category of $\TT$-modules without $\ZZ$-torsion.
Note that if $I$ and $I^{\prime}$ are ideals of $\TT$, then
$I \otimes I^{\prime}$ is isomorphic to $II^{\prime}$ under this convention, 
as the natural map sending $x \otimes x^{\prime}$ to $xx^{\prime}$ is
a surjection by construction, and the modules in question have the
same rank over $\ZZ$.  Thus the kernel is finite, and hence trivial
since under our convention $I \otimes I^{\prime}$ is torsion free.

We will find the following observation about these two operations to be
useful in what follows:
\begin{lemma} \label{lemma:star}
Let $M$ and $N$ be finitely generated, torsion free $\TT$-modules. 
We have natural isomorphisms:
\begin{enumerate}
\item $M^* \cong \Hom_{\TT}(M, \TT^*)$
\item $(M \otimes N)^* \cong \Hom_{\TT}(M,N^*)$
\item $\Hom_{\TT}(M,N)^* \cong M \otimes N^*.$
\end{enumerate}
\end{lemma}
\proof 
\begin{enumerate}
\item We have $M^* \cong \Hom_{\TT}(\TT, M^*)$.  Taking $\ZZ$-duals,
we find that
the latter is naturally isomorphic to $\Hom_{\TT}(M, \TT^*)$.  
\item By (1), 
$(M \otimes N)^*$ is canonically isomorphic to 
$\Hom_{\TT}(M \otimes N, \TT^*),$ which in turn 
is canonically isomorphic to $\Hom_{\TT}(M, \Hom_{\TT}(N, \TT^*))$,
and the result follows.  
\item This follows from (2) by replacing $N$ with $N^*$ and dualizing.\endproof
\end{enumerate}
\section{Morphisms ``up to support on $S$''} \label{sec:invert}

In what follows, we will often find it convenient to work in 
categories of abelian varieties in which certain morphisms
have been formally inverted.  In particular, we will fix a finite flat
$\ZZ$-algebra $\TT$ and an abelian variety $J_{/\QQ}$ on which
$\TT$ acts faithfully. 

\begin{defn} 
Let $\CA$ be the category whose objects consist of abelian
varieties $J^{\prime}$ with an action of $\TT$, such that
there exists a $\TT$-equivariant isogeny
$\phi: J \rightarrow J^{\prime}$ defined over $\QQ$, and whose morphisms
are $\TT$-equivariant maps of abelian varieties defined over $\QQ$.
\end{defn}

We are interested in the behavior of maps in this category up to support on 
a certain finite set $S$ of maximal ideals of $\TT$.  In order to exclude
phenomena that happen solely at maximal ideals inside $S$, we wish
to formally invert morphisms in $\CA$ which are isomorphisms at all
$\mm$ outside $S$.  

This motivates us to consider the category $\CA_S$ of ``abelian varieties
up to S-isogeny'', whose objects
are the same as the objects of $\CA$, but whose morphisms include
formal inverses for those isogenies with kernel supported on $S$.
In order to make this precise, we invoke the following lemma:

\begin{lemma} \label{lemma:TS}
Let $S$ be a finite set of maximal ideals of $\TT$.  Then
there is an element $\sigma_S$ of $\TT$ such that for any maximal
ideal $\mm$ of $\TT$, $\sigma_S \in \mm$ if and only if $\mm \in S$.
In particular, $U = \spec \TT - S$ is affine, and equal to
$\spec \TT[\frac{1}{\sigma_S}]$.
\end{lemma}

\proof
It suffices to show this when $\TT$ is reduced, as if
$\sigma_S^{\red}$ is an element of $\TT^{\red}$ (the reduced quotient
of $\TT$) with the desired property, then any lift of $\sigma_S^{\red}$
to $\TT$ will have the desired property.

Suppose $\TT$ is reduced, and let $\hat \TT$ be its normalization.
Then $\hat \TT$ is a product of maximal orders in number fields,
and so its Picard group is finite.  In particular, for each
maximal ideal $\mm$ of $\TT$, we can find an element
$\hat \sigma_{\mm}$ of $\hat \TT$ such that $\hat \sigma_{\mm}$
lies in $\mm \hat \TT$ but does not lie in $\mm^{\prime} \hat \TT$
for any $\mm^{\prime}$ different from $\mm$.

\begin{lemma} For $n$ sufficiently large, $\hat \sigma_{\mm}^n$ lies
in $\TT$.
\end{lemma}
\proof
Let $N$ be the index of $\TT$ in $\hat \TT$.  It suffices to check
this statement locally at each maximal ideal of $\TT$ containing $N$.
Since $\hat \sigma_{\mm}$ lies in $\mm \hat \TT$, for some $n$
$\sigma_{\mm}^n$ lies in $N\hat \TT_{\mm}$ and therefore lies in
$\TT_{\mm}$.  Fix a maximal ideal $\mm^{\prime}$ of $\TT$ distinct from
$\mm$ but containing $N$.  Then $\hat \sigma_{\mm}$ lies in
$\TT_{\mm^{\prime}}^{\times}$.  Let $n$ be the order of
$(\TT_{\mm^{\prime}}/N\TT_{\mm^{\prime}})^{\times}$.  Then
$\hat \sigma_{\mm}^n \in 1 + N\TT_{\mm^{\prime}} \subset \TT_{\mm^{\prime}},$
as required.  
\endproof

We now resume the proof of Lemma~\ref{lemma:TS}.  Let
$\sigma_{\mm} = \hat \sigma_{\mm}^n$.  Then $\sigma_{\mm}$ lies
in $\mm$, but does not lie in $\mm^{\prime}$ for any $\mm^{\prime}$
distinct from $\mm$.  Taking $\sigma_S$ equal to the product of
the $\sigma_{\mm}$ gives us a $\sigma_S$ with the desired properties.
\endproof

Let $\TT_S = \TT[\frac{1}{\sigma_S}]$.  We define $\CA_S$ to be
the category whose objects are the objects of $\CA$, but whose
morphisms are given by $\Hom_{\CA_S}(A_1,A_2) = 
\Hom_{\CA}(A_1,A_2) \otimes_{\TT} \TT_S$.  (The $\TT$-module structure
on $\Hom_{\CA}(A_1,A_2)$ is defined by composition either on the
left or the right; since morphisms in $\CA$ are $\TT$-equivariant it
does not matter which.) 

This has the desired effect of ignoring behavior supported on $S$,
as any morphism whose kernel is supported on $S$ factors through
a sufficiently high power of $\sigma_S$, and hence is an isomorphism
in $\CA_S$.  We will make use of this construction in the next section.

Observe that since $\sigma_S$ is a unit in $\TT_{\mm}$ for any $\mm$
outside $S$, any morphism $A_1 \rightarrow A_2$ in $\CA_S$ induces
a well defined map $A_1[\mm^{\infty}] \rightarrow A_2[\mm^{\infty}]$
which sends $x \in A_1[\mm^{\infty}]$
to $\sigma_S^{-n} (\sigma_S^nf)(x)$, for $n$ sufficiently large to
make $\sigma_S^n$ a morphism in $\CA$.

It will be convenient to speak of the ``kernel'' of a map in $\CA_S$,
even though strictly speaking this notion is not well defined.  It is
possible to make it precise ``up to support on $S$'', however;
for any morphism $f: A_1 \rightarrow A_2$ in $\CA_S$, we can fix an $n$
such that $\sigma_S^n f$ is a morphism in $\CA$, and take 
$\ker_S f = \ker \sigma_S^n f$.  This is independent of $n$ up to
a finite group supported on $S$. 

\begin{defn} Let $B$ be a subvariety of $A_1$ which is
stable under $\TT$ and closed under addition and
inverses, and $f: A_1 \rightarrow A_2$ a morphism in
$\CA_S$.  Then $\ker_S f = B$ ``up to support on $S$'' if
for some $n$, $\sigma_S^n f$ is a morphism in $\CA$, $\ker
\sigma_S^n f$ contains $B$, and the group $(\ker \sigma_S^nf)/B$
is finite and supported on $S$.  (Equivalently, $\ker_S f = B$
``up to support on $S$'' if $B[\mm^{\infty}]= \ker f: A_1[\mm^{\infty}]
\rightarrow A_2[\mm^{\infty}]$ for all $\mm$ outside $S$. ) 
\end{defn}

\section{The isogeny class of an abelian variety with $\TT$-action}
\label{sec:isogeny}

We now begin the program described in the introduction, in which we
study the category of abelian varieties isogenous to $J^D(\Gamma)$ for
some $D$ and $\Gamma$.  Our techniques apply in considerably more generality,
however, so rather than working explicitly with $J^D(\Gamma)$, 
we work in the following setting:
let $J$ be an abelian variety, defined over $\QQ$, and suppose that
$J$ comes equipped with the action of a finite, flat $\ZZ$-algebra $\TT$,
in which the action of any element $\sigma$ of $\TT$ on $J$ is defined
over $\QQ$.

If $\mm$ is a maximal ideal of $\TT$ of residue characteristic $l$, 
let $T_{\mm}J$ denote the $\mm$-adic {\em contravariant} Tate module of
$J$ at $\mm$; that is, the inverse limit of 
$J[\mm^k]^{\vee}$ over $k$, where $(-)^{\vee}$ denotes Cartier duality.  
Thus $T_{\mm}J$ is naturally isomorphic to 
$\Hom(J[\mm^{\infty}], \mu_{l^{\infty}})$.  If $\alpha \in T_{\mm}J$
and $x \in J[\mm^{\infty}]$, we denote by $\alpha(x)$ the root of
unity obtained by evaluating the element of 
$\Hom(J[\mm^{\infty}], \mu_{l^{\infty}})$ which corresponds to $\alpha$ 
at the element $x$.  For a morphism $f:A_1\rightarrow A_2$ of
abelian varieties with $\TT$-action, we let $f_{\mm}$ denote the
induced map $T_{\mm} A_2 \rightarrow T_{\mm} A_1$.

The action of $\TT$ on $J$ induces an action of $\TT \otimes \QQ$ on the 
singular cohomology $H^1(J(\CC), \QQ)$ of $J$; we assume that this makes
$H^1(J(\CC), \QQ)$ into a free $\TT \otimes \QQ$-module, of rank two.
In particular, this means that the action of $\TT$ on $J$ is {\em faithful}.

Since there are natural isomorphisms $H^1(J(\CC), \ZZ)_{\mm} \cong T_{\mm} J$,
$T_{\mm} J \otimes \QQ$ is free of rank two as a 
$\TT_{\mm} \otimes \QQ$-module. 
The action of $G_{\QQ}$ on $T_{\mm}J \otimes \QQ$ gives us a two-dimensional
representation $\rho_{\mm}$ of $G_{\QQ}$ over $\TT_{\mm} \otimes \QQ$.
We impose the further hypothesis that $\tr \rho_{\mm}(g)$ and 
$\det \rho_{\mm}(g)$
(a priori elements of $\TT_{\mm} \otimes \QQ$), lie in $\TT_{\mm}$
for all $\mm$ and $g$.  We further assume that the 
subalgebra of $\TT_{\mm}$ generated by the elements $\tr \rho_{\mm}(g)$
and $\det \rho_{\mm}(g)$ as $g$ varies is {\em reduced}
for all $\mm$.

Finally, we assume that there exists a $\TT$-equivariant map
$\psi: J \rightarrow J^{\vee},$ defined over a finite extension $K$
of $\QQ$, and a character $\chi$ of $\gal(K/\QQ)$,
taking values in $\TT^{\times}$,
such that $\psi$ has finite kernel, and 
for all $\sigma \in \gal(K/\QQ)$, we have
$\sigma \psi= \chi(\sigma) \psi \sigma$.

\begin{remark} \rm \label{remark:modular} 
The motivation for these assumptions comes from the
example in which $J=J^D(\Gamma)$ for some discriminant $D$ and
congruence subgroup $\Gamma$ of level $N$, and $\TT$ is the Hecke algebra 
$\TT^D(\Gamma).$ In this case it can be shown 
(see for instance Carayol~\cite{Ca})
that $H^1(J(\CC),\QQ)$ is free of rank two over $\TT \otimes \QQ$.
Moreover, for a prime $l$ not dividing $ND$, $\tr \rho_{\mm}(\Frob_l)$
is equal to $T_l \in \TT$, and $\det \rho_{\mm}(\Frob_l) = l\<l\>$.
The subalgebra of $\TT$ generated by $\<l\>$ and
$T_l$ for $l$ prime to $ND$ is indeed reduced, as it is well-known that
there exists a basis of simultaneous eigenforms for all such $T_l$.
Since the elements $\Frob_l$ for $l$ prime to $N$ are dense in $G_{\QQ}$,
and $\TT_{\mm}$ is complete, the subalgebra of $\TT_{\mm}$ generated
by the traces and determinants of elements of $G_{\QQ}$ is generated
by the aforementioned $\<l\>$ and $T_l$, and is therefore reduced.
Finally, the map $\psi$ described above is constructed in 
Lemma~\ref{lemma:rosati}.
It follows that the hypotheses of this section are satisfied in this case.
\end{remark}

The existence of the map $\psi$ described above has strong consequences
for the structure of $\TT$.  In particular, $\psi$ induces an
isomorphism $T_{\mm} J \otimes \QQ \rightarrow (T_{\mm} J)^* \otimes \QQ$.
Since the former is free of rank two over $\TT_{\mm} \otimes \QQ$,
and the latter is isomorphic to $(\TT_{\mm}^* \otimes \QQ)^2$, it follows
that $\TT_{\mm}^* \otimes \QQ$ is free of rank one over
$\TT_{\mm} \otimes \QQ$.  In particular, $\TT \otimes \QQ$ is a 
zero-dimensional Gorenstein ring.

Following the analysis of the Tate modules of modular Jacobians 
in~\cite{Ca}, we
fix a maximal ideal $\mm$ of $\TT$, and let $\TT^{\prime}_{\mm}$ 
denote the subalgebra of $\TT_{\mm}$ generated by the traces and
determinants of
$\rho_{\mm}(g)$ for $g \in G_{\QQ}$.  By our hypothesis, this is
a reduced subalgebra of $\TT_{\mm}$.  We can thus consider the
normalization $\hat \TT^{\prime}_{\mm}$ of $\TT^{\prime}_{\mm}$;
this will be a product of discrete valuation rings $\OO_i$, with
corresponding fields of fractions $K_i$.

For each $i$, $T_{\mm} J \otimes_{\TT_{\mm}^{\prime}} K_i$ is
free of rank two over the Artinian $K_i$-algebra 
$\TT_{\mm} \otimes_{\TT_{\mm}^{\prime}} K_i$; we consider the
representation $\rho_{\mm,i}$ of $G_{\QQ}$ over 
$\TT_{\mm} \otimes_{\TT_{\mm}^{\prime}} K_i$ obtained in this manner. 
Since the traces and determinants
of $\rho_{\mm}$ lie in $\TT_{\mm}^{\prime}$, the
traces and determinants of $\rho_{\mm,i}$ lie in $\OO_i$ for all $i$.  
Thus $\rho_{\mm,i}$
can be considered as a representation of $G_{\QQ}$ over $\OO_i$.

Let $k_i$ denote the residue field of $\OO_i$; then we may consider the
mod $\mm$ reduction $\overline{\rho}_{\mm,i}: G_{\QQ} \rightarrow
\GL_2(k_i)$ of $\rho_{\mm,i}$.  
For $g \in G$, $\tr \overline{\rho}_{\mm,i}(g)$ is
equal to the reduction of $\tr \rho_{\mm}(g)$ modulo $\mm$, and
is therefore independent of $i$.  Thus the semisimplification of
$\overline{\rho}_{\mm,i}$ is independent of $i$. 

With this definition, we make one final assumption on $J$, namely,
that for all but finitely many $\mm$, the representations
$\overline{\rho}_{\mm,i}$ are absolutely irreducible.
For such $\mm$, $\overline{\rho}_{\mm,i}$
is independent of $i$ and we denote it simply by $\overline{\rho}_{\mm}$.

\begin{remark} \rm It should be pointed out here that in the case in
which $J$ is a modular Jacobian, the $\overline{\rho}_{\mm}$ constructed 
above is the Cartier dual of the representation usually called 
$\overline{\rho}_{\mm}$ in the literature.  (This is because the latter
is defined in terms of the $\mm$-torsion of $J$, whereas the
representations we consider are defined in terms of $T_{\mm}J$ instead.)
\end{remark}

Fix a finite set $S$ of maximal ideals of $\TT$, containing all $\mm$
for which one (and hence all) of the $\overline{\rho}_{\mm,i}$ is not absolutely
irreducible.  The following result is a straightforward generalization 
of~\cite{Ca}, Theorem 3:

\begin{prop} \label{prop:ca} 
For all $\mm$ outside $S$, there exists a $G_{\QQ}$-module
$V_{\mm}$, free of rank two over $\TT_{\mm}$, such that $V_{\mm} \otimes \QQ$
and $T_{\mm} J \otimes \QQ$ are isomorphic as $G_{\QQ}$-modules.
\end{prop}
\proof Consider the action of $G_{\QQ}$ on 
$\hat{V}_{\mm}^{\prime} = (\hat{\TT}^{\prime}_{\mm})^2$
in which $G_{\QQ}$ acts on each factor $(\OO_i)^2$ of 
$(\hat{\TT}^{\prime}_{\mm})^2$ by $\rho_{\mm,i}$.  This
gives a representation of $G_{\QQ}$ in $\GL_2(\hat{\TT}^{\prime}_{\mm})$
whose character agrees with $\rho_{\mm}$.  
By~\cite{Ca}, Theorem 2, this
representation is realizable over $\TT_{\mm}^{\prime}$; i.e., there is
a $G_{\QQ}$-module $V_{\mm}^{\prime}$ free of rank two over 
$\TT^{\prime}_{\mm}$ such that the character of the $G_{\QQ}$-action agrees
with that of $\rho_{\mm}$.  Then 
$V_{\mm} = V_{\mm}^{\prime} \otimes_{\TT_{\mm}^{\prime}} \TT_{\mm}$ is free
of rank two over $\TT_{\mm}$, and $V_{\mm} \otimes \QQ$ is equivalent to
$\rho_{\mm}$ by~\cite{Ca}, Theorem 1.
\endproof 

Our goal is to study $\TT$-equivariant isogenies of $J$ which are
defined over $\QQ$, up to support on $S$.  The first step in this study
is to understand the structure of $T_{\mm}J$ for $\mm$ outside $S$.

\begin{lemma} \label{lemma:multiplicity} 
For all but finitely many maximal ideals $\mm$ of $\TT$,
$T_{\mm}J$ is free of rank two over $\TT_{\mm}$.
\end{lemma}
\proof
Since we have assumed that $H^1(J(\CC),\QQ)$ is free of rank two over
$\TT \otimes \QQ$, we can choose two elements of $H^1(J(\CC),\ZZ)$
that generate $H^1(J(\CC),\QQ)$ over $\TT \otimes \QQ$.  These two
elements give us a map $\TT^2 \rightarrow H^1(J(\CC), \ZZ)$, with
finite cokernel.  Localizing this map at any $\mm$ on which the cokernel
has no support gives the desired isomorphism $\TT_{\mm}^2 \cong T_{\mm}J$.
\endproof

We will see (Lemma~\ref{lemma:submodule}, below) that the Tate modules
at such $\mm$ are particularly well-behaved.  It will thus be convenient
for us to replace $J$ with an isogenous abelian variety, which we call
$\Jmin$, for which $T_{\mm} \Jmin$ is free of rank two over $\TT_{\mm}$
for all $\mm$ outside $S$.  The following lemma allows us to do so.
 
\begin{lemma} \label{lemma:Jmin}
There exists an abelian variety $\Jmin$, with a $\TT$-action satisfying
the above hypotheses, such that
$J$ and $\Jmin$ are $\TT$-equivariantly isogenous 
(over $\QQ$), and $T_{\mm}\Jmin$ is free of rank two over
$\TT_{\mm}$ for all $\mm$ outside $S$. 
\end{lemma}
\proof
Let $S^{\prime}$ be the set of $\mm$ outside $S$ at which $T_{\mm} J$
is not free of rank two; then $S^{\prime}$ is finite.  For
each $\mm$ in $S$, fix an isomorphism 
$V_{\mm} \otimes \QQ \cong T_{\mm} J \otimes \QQ$.  Multiplying this
isomorphism by a sufficiently large $n \in \ZZ$, we may assume that it
maps $V_{\mm}$ into $T_{\mm} J$, with finite index.  Let 
$W_{\mm} \subset J[\mm^{\infty}]$ be the finite subset of the $\mm$-divisible
group annihilated by the image of $V_{\mm}$ in $T_{\mm} J$, under the natural
pairing
$$T_{\mm} J \times J[\mm^{\infty}] \rightarrow \mu_{p^{\infty}}.$$  If we take
$\Jmin$ to be the quotient of $J$ by the sum of the $W_{\mm}$, we find
that the map $T_{\mm} \Jmin \rightarrow T_{\mm} J$ identifies $T_{\mm} \Jmin$
with the image of $V_{\mm}$ for all $\mm$ in $S^{\prime}$, and is an 
isomorphism for $\mm$ outside $S^{\prime}$; in either case $T_{\mm} \Jmin$
is free of rank two over $\TT_{\mm}$. 
\endproof

We fix, for the remainder of this section, a $\Jmin$ as above.

\begin{lemma} \label{lemma:submodule}
Let $M$ be a $G_{\QQ}$-stable submodule of $T_{\mm} \Jmin$, of finite
index.  If
$\mm$ lies outside $S$, then $M = IT_{\mm}\Jmin$ for some ideal $I$ of $\TT$.
\end{lemma}
\proof
We use induction on the length of a maximal $G_{\QQ}$-stable filtration
of $T_{\mm}\Jmin/M$.  If this length is zero, then $I$ is the unit ideal and
the result is clear.  In general, let 
$M = M_n \subset M_{n-1} \subset \dots \subset M_0 = T_{\mm}\Jmin$
be a maximal $G_{\QQ}$-stable chain of submodules.  By induction we may
assume that $M_{n-1} = I^{\prime}T_{\mm}\Jmin$ for some $I^{\prime}$.  

Consider $\mm M_{n-1} + M$.  This is $G_{\QQ}$-stable, contains $M$,
and is contained in $M_{n-1}$.  It is thus equal to either $M$ or
$M_{n-1}$.  Suppose it were equal to $M_{n-1}.$  Then $M/\mm M_{n-1}$
would equal $M_{n-1}/\mm M_{n-1}$, which is impossible since then
$M$ would contain a set of generators for $M_{n-1}$ by Nakayama's lemma.
Thus $\mm M_{n-1} + M = M$, so $M$ contains $\mm M_{n-1}$.

Let $V_{\overline{\rho}_{\mm}}$ be a two-dimensional $(\TT/\mm)$
vector space on which $G_{\QQ}$ acts via $\overline{\rho}_{\mm}$.
The module $M_{n-1}/\mm M_{n-1}$ is $G_{\QQ}$-equivariantly isomorphic to 
$(I^{\prime}/\mm I^{\prime}) \otimes_{\TT/\mm} V_{\overline{\rho}_{\mm}},$
where $G_{\QQ}$ acts trivially on $I^{\prime}/\mm I^{\prime}$.  Let $V$ be
the image of $M$ in $M_{n-1}/\mm M_{n-1}$.  Since $V$ is $G_{\QQ}$-invariant,
and $V_{\overline{\rho}_{\mm}}$ is irreducible, $V$ is given by
$\hat{V} \otimes V_{\overline{\rho}_{\mm}}$ for some subspace
$\hat{V}$ of $I^{\prime}/\mm I^{\prime}$.  Let $I$ be the preimage of
$\hat{V}$ in $I^{\prime}$.  Then $IT_{\mm}J = M$, since both
contain $\mm M_{n-1}$ and map to $V$ modulo $\mm M_{n-1}$.
\endproof

This result has important consequences for $\TT$-equivariant
maps from $\Jmin$.  In particular, for any abelian varieties $A_1, A_2$ which
are $\TT$-equivariantly isogenous to $J$ over $\QQ$, let $\Hom(A_1,A_2)$
(resp. $\End(A_1)$) denote the $\TT$-module of $\TT$-equivariant
maps: $A_1 \rightarrow A_2$ (resp. $\TT$-equivariant endomorphisms of $A_1$)
which are defined over $\QQ$.  Then:

\begin{cor} \label{cor:end}
The cokernel of the natural map $\TT \rightarrow \End(\Jmin)$ is
supported on $S$.
\end{cor}
\proof Fix $\mm$ outside $S$; we show that the natural
map 
$$\TT_{\mm} \rightarrow \End(\Jmin)_{\mm}$$ is an isomorphism.
We have an isomorphism $$\End(\Jmin)_{\mm} \cong \End_G(T_{\mm} \Jmin),$$
by~\cite{Fa}, Theorem 4.
$\End_G(T_{\mm} \Jmin)$ is clearly faithful over $\TT_{\mm}$, and finitely
generated.  It thus suffices to show that 
$$\End_G(T_{\mm} \Jmin) \otimes \TT/\mm$$ is 1-dimensional.  
Observe that $\End_G(\Jmin[\mm])$ is 1-dimensional, by Schur's lemma and
the (absolute) irreducibility of $\overline{\rho}_{\mm}$.  We construct 
an injection: 
$$\End_G(T_{\mm} \Jmin) \otimes \TT/\mm \rightarrow \End_G(\Jmin[\mm]),$$
the existence of which immediately implies the desired result.
The exact sequence
$$0 \rightarrow \mm T_{\mm}\Jmin \rightarrow T_{\mm}\Jmin \rightarrow
\Jmin[\mm]^{\vee} \rightarrow 0$$
induces an exact sequence:
$$0 \rightarrow \Hom_G(T_{\mm}\Jmin, \mm T_{\mm}) \rightarrow 
\End_G(T_{\mm}\Jmin) \rightarrow 
\Hom_G(T_{\mm}\Jmin, \Jmin[\mm]^{\vee}),$$ 
and the latter is isomorphic to $\End_G(\Jmin[\mm])$ via Cartier duality.
Now 
$\Hom_G(T_{\mm}\Jmin, \mm T_{\mm})$ clearly contains
$\mm\End_G(T_{\mm} \Jmin)$, and it suffices to show that they are equal.

Let $\phi \in \Hom_G(T_{\mm} \Jmin, \mm T_{\mm})$; by 
Lemma~\ref{lemma:submodule}, the image of $\phi$ is isomorphic to
$IT_{\mm}\Jmin$ for some $I \subset \mm$.  Since 
$IT_{\mm} \Jmin \cong T_{\mm} \Jmin$ and $T_{\mm} \Jmin$ is free over 
$\TT_{\mm}$, $I$ is principal; say generated by $\sigma$.  Thus $\phi/\sigma$
is a well defined endomorphism of $T_{\mm} \Jmin$.  Since $\sigma \in \mm$,
$\phi$ lies in $\mm\End_G(T_{\mm} \Jmin)$ as required.
\endproof
 
\begin{cor} \label{cor:kernel} 
Let $\phi: \Jmin \rightarrow A$ be a $\TT$-equivariant
isogeny.  Let $I = \ann(\ker \phi)$.  Then 
$(\ker \phi)_{\mm} = \Jmin[I]_{\mm}$ for all $\mm$ outside $S$.
\end{cor}
\proof
The map $\phi$ induces an exact sequence:
$$0 \rightarrow T_{\mm} A \rightarrow T_{\mm} \Jmin \rightarrow
    (\ker \phi)_{\mm}^{\vee} \rightarrow 0,$$
In particular the image of 
$T_{\mm} A$ is a $\TT$-stable submodule of $T_{\mm} \Jmin$.  
By Lemma~\ref{lemma:submodule}, we can find an ideal $I^{\prime}$ of $\TT$
such that the image of $T_{\mm} A$ is 
$I^{\prime}_{\mm}T_{\mm}\Jmin$ for all $\mm$ outside $S$.

Now $[\ker \phi]_{\mm}^{\vee} = T_{\mm}\Jmin/I^{\prime}T_{\mm}\Jmin$ for all 
$\mm$ outside $S$, by the above exact sequence.  Thus 
$I^{\prime}_{\mm} = I_{\mm}$
for all $\mm$ outside $S$.  Since $T_{\mm}\Jmin/IT_{\mm}\Jmin =
\Jmin[I]_{\mm}^{\vee}$, the result follows. 
\endproof

Note that the above result implies that the isomorphism class of $I_{\mm}$
does not depend on the isogeny $\phi$, but only on the variety $A$.
In fact, this is true globally, as well as locally, as long as one works
``up to support on $S$''.  To see this, fix a 
$\phi$ and $I$ as above.  Then
for $\sigma \in I$, $\ker \phi \subset \Jmin[I] \subset \Jmin[\sigma]$.  
Thus the map
$\sigma: \Jmin \rightarrow \Jmin$ factors through $\phi$; i.e. 
$\sigma = \varphi \phi$
for a unique map $\varphi: A \rightarrow \Jmin$.  Let $f$ to be the
map which associates to each $\sigma$ the corresponding map $\varphi$.

\begin{prop} \label{prop:functor1}
The map $f: I \rightarrow \Hom(A, \Jmin)$ is injective, 
and $f_{\mm}$ is an isomorphism for all $\mm$ outside $S$.
\end{prop}
\proof
The map $f$ must be injective, as if $\sigma$ induced the zero map
$A \rightarrow \Jmin$, then the 
action of $\TT$ on $\Jmin$ could not be faithful.  
It remains to show that $f_{\mm}$ is 
surjective for $\mm$ outside $S$.

Let $\varphi$ be an element of $\Hom(A,\Jmin)$.  Then the composition
$\varphi \phi$ lies in $\Hom(\Jmin,\Jmin)$.  
By Corollary~\ref{cor:end}, the inclusion 
$\TT_{\mm} \rightarrow \Hom(\Jmin,\Jmin)_{\mm}$ is an isomorphism.  Thus we 
can choose
an element $\sigma$ of $\TT_{\mm}$ corresponding to $\varphi\phi$; since 
$\Jmin[I]_{\mm}$ 
is contained in the kernel of $\varphi\phi$ by Corollary~\ref{cor:kernel}, 
$\sigma$ lies in $I_{\mm}$.  Moreover, $f_{\mm}(\sigma) = \varphi$ by 
construction.
\endproof

\begin{cor} \label{cor:tate} 
If $A$ is $\TT$-equivariantly isogenous to $J$,
and $\mm$ lies outside $S$, 
then the natural map
$T_{\mm} \Jmin \otimes \Hom(A,\Jmin)_{\mm} \rightarrow
T_{\mm} A$ is an isomorphism.  In particular,
$T_{\mm} A$ is isomorphic to $\Hom(J^{\prime},\Jmin)_{\mm}^2$ as a
$\TT_{\mm}$-module.  
\end{cor}
\proof Since $T_{\mm} \Jmin$ is free over $\TT_{\mm}$,
$T_{\mm} \Jmin \otimes \Hom(A,\Jmin)_{\mm}$ is torsion free.
It thus suffices to show that the natural map is surjective, as the
modules in question have the same $\ZZ$-rank.

Fix an isogeny $\phi: \Jmin \rightarrow A$.  Then
$\phi$ identifies $T_{\mm}A$ with $IT_{\mm}\Jmin$, where
$I = \ann (\ker \phi)$.  Let $x \in T_{\mm}A$; then
$\phi$ identifies $x$ with an element $y$ of $IT_{\mm}\Jmin$;
we write $y$ as a sum of elements of the form
$\sigma_i y_i$ with $\sigma_i \in I$ and $y_i \in T_{\mm}\Jmin$.
Since $\sigma_i \in I$, the map $\sigma_i: \Jmin \rightarrow \Jmin$
factors through $\phi$, i.e., we have maps $\varphi_i$ such that
$\sigma_i = \varphi_i\phi$.  Then $\sigma_i y_i = \varphi_i(\phi(y_i))$,
so $\phi(y_i) \otimes \varphi_i$ maps to $\sigma_i y_i$ under the natural
map.  Thus each of the $\sigma_i y_i$ are in the image, so $y$ is
in the image as well.  The second statement follows immediately since
$T_{\mm} \Jmin$ is free of rank two over $\TT_{\mm}$.
\endproof

Thus to any $A$ which is $\TT$-equivariantly isogenous to $J$ we
can associate the module $\Hom(A,\Jmin)$.  This module is
torsion free as a $\ZZ$-module, 
and $\Hom(A,\Jmin) \otimes \QQ$ is free of rank one
as a $\TT \otimes \QQ$-module by Proposition~\ref{prop:functor1}.  
Following the previous section, we
let $\CA$ be the category whose objects consist of abelian
varieties $A$ such that there exists a $\TT$-equivariant isogeny
$\phi: J \rightarrow A$ defined over $\QQ$, and whose morphisms
are $\TT$-equivariant maps of abelian varieties defined over $\QQ$.  
Then $A \mapsto \Hom(A, \Jmin)$ is a well defined
functor $\CA \rightarrow \CM$, where
$\CM$ denotes the category of finitely generated $\TT$-modules $M$ such that 
$M$ is torsion free
over $\ZZ$ and $M \otimes \QQ$ is free of rank one over $\TT \otimes \QQ$.

Let $\TT_S$ be the $\TT$-algebra defined in the previous section,
and consider the category $\CA_S$ defined therein.
The functor defined above induces a functor $\CA_S \rightarrow \CM_S$,
where $\CM_S$ is the category of ``rank one'' $\TT_S$-modules,
i.e., the category of finitely generated $\TT_S$-modules that are
torsion free over $\ZZ$ and become free of rank one over $\TT_S \otimes \QQ$. 

\begin{thm} \label{thm:equiv}
The functor $A \mapsto [A] = \Hom(A,\Jmin)$ is 
an antiequivalence of categories $\CA_S \rightarrow \CM_S.$  
\end{thm}

\proof We must show that the functor is full, faithful, and dense.  Given
an $M \in \CM$, we can find an injection $M \rightarrow \TT$.  Let $I$
be its image.  Then by Proposition~\ref{prop:functor1}, we have a map
$I \rightarrow \Hom(\Jmin/\Jmin[I],\Jmin)$ whose cokernel is supported on $S$.
In particular, $\Hom(\Jmin/\Jmin[I],\Jmin)$ is isomorphic to $M$ in $\CM_S$.  
Thus the functor is dense.

It remains to show that the functor is fully faithful.  Since 
$\Hom_{\CA_S}(A_1,A_2)$ and $\Hom_{\CM_S}([A_2],[A_1])$
are $\TT_S$-modules, it suffices to show that the 
map $\Hom(A_1,A_2)_{\mm} \rightarrow \Hom([A_2],[A_1])_{\mm}$ is an
isomorphism for any $\mm$ outside $S$ and any $A_1,A_2 \in \CA$.
We have a commutative diagram:
$$\begin{CD}
\Hom(A_1,A_2)_{\mm} @>>> \Hom([A_2],[A_1])_{\mm} \\
@VVV @VVV\\
\Hom(T_{\mm} A_2, T_{\mm} A_1) @>>> 
\Hom(T_{\mm} \Jmin \otimes [A_2], T_{\mm} \Jmin \otimes [A_1])\\
\end{CD},$$
in which the bottom horizontal map is induced by the isomorphisms
$[A_1] \otimes T_{\mm} \Jmin \rightarrow T_{\mm} A_1$ and
$[A_2] \otimes T_{\mm} \Jmin \rightarrow T_{\mm} A_2$ of
Corollary~\ref{cor:tate}.
Note that each map appearing in this diagram is injective, as the modules
in question have no $\ZZ_l$-torsion and each map becomes an isomorphism
after tensoring with $\QQ$.  Moreover, the left-hand vertical map
and the bottom horizontal map are isomorphisms.  Thus the right-hand
vertical map is surjective, and hence an isomorphism.  It follows
that the upper horizontal map is an isomorphism as well.
\endproof

In light of this result, we will primarily be interested in objects
and morphisms in $\CA_S$ and $\CM_S$ for the remainder of this chapter.
Rather than clutter our notation by writing ``$\otimes_{\TT} \TT_S$''
repeatedly, we adopt the following convention: if $M$ and $N$ are
objects of $\CM$, an $S$-morphism: $M \rightarrow N$ is a map
of $\TT_S$-modules $M \otimes_{\TT} \TT_S \rightarrow N \otimes_{\TT} \TT_S.$
An $S$-isomorphism $M \rightarrow N$ is an $S$-morphism that is an
isomorphism in $\CM_S$, and $M$ and $N$ are said to be $S$-isomorphic
if there exists an $S$-isomorphism between them.  (In this case, we
write $M \cong_S N$.)

The formalism we have introduced gives information about isogenies in $\CA$, 
up to support on $S$, but it
does so in terms of our choice of $\Jmin$, which was somewhat arbitrary.
Our next goal is to obtain a result that does not involve this choice.
Let us assume that we have objects $A_1$ and $A_2$ in $\CA$; then we
can study the maps $A_1 \rightarrow A_2$ by studying maps from
$[A_2]$ to $[A_1]$.  In general $[A_2]$ and $[A_1]$ depend on $\Jmin$.
However, if we have a $\TT$-module $M$ such that $[A_2] \cong_S [A_1] \otimes
M$ for one choice of $\Jmin$, then it is easy to check that this relationship
is independent of the choice of $\Jmin$.  (Changing $\Jmin$ amounts to
twisting $[A_2]$ and $[A_1]$ by the same locally free $\TT$-module.) 
In this case we can characterise maps from $A_1$ to $A_2$ in terms of $M$.
Recall that we have adopted the convention that if $M$ and $M^{\prime}$
are objects of $\CM$, then $M \otimes M^{\prime}$ is taken to mean
the usual tensor product modulo $\ZZ$-torsion (and thus also lies in $\CM$).

If we have an $S$-morphism $$h: [A_1] \otimes M \rightarrow [A_2],$$ then
$h$ induces an $S$-morphism $$\hat h: M \rightarrow \Hom([A_1],[A_2]),$$ 
defined by $\hat h(m)(x) = h(m \otimes x)$ for $x \in [A_1], m \in M$.  By
Theorem~\ref{thm:equiv}, we can consider $\hat h$ to be a map
$M \rightarrow \Hom(A_2,A_1)$.  This
gives us a map $g:\Hom(A_1,A_2) \rightarrow \Hom(M,\End(A_1))$,
defined by $g(\phi)(m) = \hat h(m) \circ \phi$, for $\phi:A_1 \rightarrow A_2$. 

For an element $\phi$ of $\Hom_{\CA_S}(A_1,A_2)$,
let $I_{\phi,S} \subset \End(A_1) \otimes_{\TT} \TT_S$
be the ideal generated by the image
of $M \otimes_{\TT} \TT_S$ under $g(\phi)$, and let
$I_{\phi}$ be the ideal $I_{\phi,S}
\cap \End(A_1).$

\begin{prop} \label{prop:JJ'}
The map $g$ defined above is an $S$-isomorphism.  In particular,
$g$ induces an isomorphism $$\Hom_{\CA_S}(A_1,A_2) \rightarrow
\Hom_{\CM_S}(M,\End(A_1)).$$
Moreover, if $\phi$ is an element of $\Hom_{\CA_S}(A_1,A_2)$, then
$\ker_S \phi = A_1[I_\phi]$ ``up to support on $S$.''
\end{prop}
\proof 
We construct an inverse to $g$ in $\CM_S$.  We have a sequence of
natural $S$-morphisms: 
$$\begin{array}{lll}
\Hom(M,\End(A_1)) & \rightarrow & \Hom(M,\End([A_1])) \\ 
& \rightarrow &
\Hom(M \otimes [A_1], \End([A_1]) \otimes [A_1]) \\
& \stackrel{h}\rightarrow  
& \Hom([A_2],[A_1]) \\
& \rightarrow & \Hom(A_1,A_2),
\end{array}$$  
where the first and last morphisms come from Theorem~\ref{thm:equiv},
and the second sends a map $f$ to the map $f \otimes 1$.
It is straightforward to check that this provides an inverse to $g$. 

Now fix a $\phi \in \Hom_{\CA_S}(A_1,A_2)$, and an $\mm$ outside $S$, and
let $I_{\phi}$ be the corresponding ideal defined above. Since 
$[A_1]_{\mm} \otimes M_{\mm}$ maps onto $[A_2]_{\mm}$,
the map $T_{\mm} A_1 \otimes M_{\mm} \rightarrow T_{\mm} A_2$
defined by sending $\beta \otimes m$ to $(\hat h(m))_{\mm}(\beta)$, is
surjective.  (The latter map can be obtained from the former by
tensoring with $T_{\mm} \Jmin$ and exploiting the isomorphism
$T_{\mm} J^{\prime} \cong T_{\mm} \Jmin \otimes [J^{\prime}]$ for
$J^{\prime} = A_1$ or $A_2$.) 

Thus we can write any $\alpha \in T_{\mm} A_2$ as a sum of 
terms of
the form $$(\hat h(m_i))_{\mm}(\beta_i)$$ for some $m_i$ in $M_{\mm}$ and 
$\beta_i$ in $T_{\mm} A_1$.  But by our construction
of $g$, we have $$\phi_{\mm}(\hat h(m_i))_{\mm}(\beta_i) 
= g(\phi)(m_i)\beta_i,$$ and the latter 
is in $I_{\phi} T_{\mm} A_1$ for all $i$.  Thus 
$\phi_{\mm}$ maps $\alpha$ into $I_{\phi} T_{\mm} A_1$.  Hence
the image of $\phi_{\mm}$ is contained in $I_{\phi,\mm} T_{\mm} A_1$.

Conversely, observe that by construction, every element of $I_{\phi}$
factors through $\phi$, since any element of $I_{\phi}$ is in
the image of $g(\phi)$ and the definition of $g(\phi)$ involved
composition with $\phi$.  Thus the image of $\phi_{\mm}$ contains
$I_{\phi} T_{\mm} A_1$.

It follows that $\phi_{\mm} T_{\mm} A_2 = I_{\phi} T_{\mm} A_1$
for all $\mm$ outside $S$.  Since the former is the annihilator
of $\ker \phi: A_1[\mm^{\infty}] \rightarrow A_2[\mm^{\infty}]$ and the 
latter is the annihilator of $A_1[\mm^{\infty}][I_{\phi}]$, we have 
$\ker_S \phi = A_1[I_{\phi}]$ up to support on $S$.  \endproof 

The behavior of certain invariants of $J$ under isogeny can be studied
using this formalism.  In particular,
let $p$ be a prime, and suppose that $J$ has purely toric reduction at $p$.
Then the special fiber of the Ner\'on model of $J$ at $p$ is an extension 
of a finite group scheme by a torus $T$, and we can consider the character 
group $X_p(J) = \Hom_{\overline{\FF}_p}(T,\GG_m)$ of $T$.  
More generally, any object $A$ of $\CA$ will also
have purely toric reduction at $p$, and we can consider its character group
$X_p(A)$.

\begin{lemma} The character group $X_p(A)$ is a ``rank one''
$\TT$-module. 
\end{lemma}
\proof
By construction, $X_p(A)$ has finite rank over $\ZZ$, and is torsion free
as a $\ZZ$-module.  It thus suffices to show that $X_p(A) \otimes \QQ$
is free of rank one as a $\TT \otimes \QQ$-module.  Let $\pp$ be a maximal
ideal of $\TT \otimes \QQ$.  It is then enough to show that
$X_p(A)_{\pp}$ is free of rank one over $\TT \otimes \QQ$.

Since $\TT \otimes \QQ$ is Gorenstein, $(\TT \otimes \QQ)[\pp]$ is
generated (as a $(\TT \otimes \QQ)/\pp$-vector space) by a
single element $\sigma$.  Since $\TT$ acts faithfully on $J$, it
acts faithfully on $A$ as well, and hence on $X_p(A)$.  In particular
there is an $x$ in $X_p(A)_{\pp}$ not killed by $\sigma$.  Let $I$ be
the annihilator of $x$ in $\TT \otimes \QQ$.  Then $I$ is killed by some
power of $\pp$, so $I[\pp]$ is nontrivial unless $I$ is the zero ideal.
But if $I[\pp]$ were nontrivial, then $I$ would contain $\sigma$, which
does not annihilate $x$ by assumption.  Thus $I$ is the zero ideal, and
$X_p(A)_{\pp}$ contains a free $(\TT \otimes \QQ)_{\pp}$-submodule
for all $\pp$.

Counting $\QQ$-ranks, we find that $X_p(A)_{\pp}$ is free of rank
one over $(\TT \otimes \QQ)_{\pp}$ for all $\pp$, since the
rank of $X_p(A)$ is equal to the dimension of $A$, and hence to
the $\ZZ$-rank of $\TT$.
\endproof

If $A_1$ and $A_2$ are two objects in $\CA$, and
$\phi: A_1 \rightarrow A_2$ is a $\TT$-equivariant map
defined over $\QQ$, then 
$\phi$ induces a natural map $X_p(A_2) \rightarrow X_p(A_1)$. 
By abuse of notation, we denote this map by $\phi$ as well.  (It will
always be clear from the context whether we mean a morphism of abelian
varieties or the map it induces on the character groups.)

\begin{prop} \label{prop:charisogeny}
Let $A \in \CA$.  The natural map 
$$X_p(\Jmin) \otimes \Hom(A,\Jmin) \rightarrow X_p(A)$$
defined by $\phi \otimes x \mapsto \phi(x)$ is an isomorphism in $\CM_S$. 
\end{prop}
\proof
Let $T_l A$ denote the {\em covariant} $l$-adic Tate module
of $A$, and $(T_l A)^t$ be the submodule of $T_l A$ defined
in~\cite{SGA7}, chapter 5.  By Grothendieck's orthogonality theorem, 
$(T_l A)^t$ is the orthogonal complement of a certain submodule
$(T_l A^{\vee})^f$ of $T_l A^{\vee}$ under the Weil pairing;
in particular if $l^nx \in (T_l A)^t$ for some x in $T_l A$
and $n \in \ZZ$, then $x \in (T_l A)^t$ as well.  It follows that
the map $\Hom(T_l A, \ZZ_l) \rightarrow \Hom((T_l A)^t, \ZZ_l)$
is surjective.  Localizing at $\mm$, we obtain a surjection
$T_{\mm} A \rightarrow \Hom((T_l A)^t, \ZZ_l)_{\mm}$, where as
usual $T_{\mm} A$ denotes the $\mm$-adic contravariant
Tate module.  Since $(T_l A)^t \cong \Hom(X_p(A), \ZZ_l(1)),$
by the results of~\cite{SGA7}, chapter 5, this gives us a functorial
surjection: $T_{\mm} A \rightarrow X_p(A)_{\mm}$.  Similarly,
we obtain a functorial surjection $T_{\mm} \Jmin \rightarrow
X_p(\Jmin)_{\mm}.$ 

These surjections fit into the commutative diagram: 
$$\begin{CD}
T_{\mm} \Jmin \otimes \Hom(A,\Jmin)_{\mm} @>>> 
T_{\mm} A \\
@VVV @VVV \\
X_p(\Jmin)_{\mm} \otimes \Hom(A,\Jmin)_{\mm} @>>> 
X_p(A)_{\mm},
\end{CD}$$
in which the horizontal maps are the natural ones, the right hand
vertical map is the surjection defined above, and the left hand vertical
map is the surjection defined above tensored with the identity.  By
Corollary~\ref{cor:tate}, the upper horizontal map is surjective;
since both vertical maps are surjections the lower horizontal map
must be surjective as well. 
\endproof

We will also have use for the following result, which is an immediate
consequence:
\begin{cor} \label{cor:controllable}
Let $\mm$ be a maximal ideal outside $S$, and suppose that 
$X_p(\Jmin)_{\mm}$ is free of rank one over $\TT_{\mm}$. 
Then for any $A \in \CA,$ 
we have $T_{\mm} A \cong (X_p(A)_{\mm})^2.$
\end{cor}
\proof 
By the proposition, and the fact that $X_p(A)_{\mm}$ is free
of rank one over $\TT_{\mm},$ we have a natural isomorphism:
$X_p(A)_{\mm} \cong [A]_{\mm}.$  Since $T_{\mm} A \cong
T_{\mm} \Jmin \otimes [A]_{\mm},$ and $T_{\mm} \Jmin$ is free
of rank two, the result follows. 
\endproof 

Finally, we wish to interpret duality of abelian varieties in
terms of the formalism we have constructed.  This presents a minor difficulty;
namely, the dual of an abelian variety $J$ in $\CA$ need not lie in $\CA$.
In fact, the abelian variety $J^{\vee}$ will be a ``twist'' of some
abelian variety in $\CA$ by the character $\chi$ we fixed at
the beginning of this section, in a sense we will make precise below.

Let $\CA_{\chi}$
denote the category of abelian varieties which are $\TT$-equivariantly
isogenous to $J^{\vee}$, and whose morphisms are $\TT$-equivariant
and defined over $\QQ$.  Duality of abelian varieties yields a 
functor $A \mapsto A^{\vee}$ from $\CA$ to $\CA_{\chi}$.

Consider also the category $\CA_{1,\chi}$, whose objects are
triples $(A,A_{\chi}, \varphi)$, where $A \in \CA$, $A_{\chi} \in \CA_{\chi}$,
and $\psi$ is a $\TT$-equivariant isomorphism $A \rightarrow A_{\chi}$,
defined over $K$, such that $\sigma \varphi = \chi(\sigma) \varphi \sigma$,
for all $\sigma \in \gal(K/\QQ)$.  A morphism in $\CA_{1,\chi}$ from
$(A,A_{\chi}, \varphi)$ to $(A^{\prime},A^{\prime}_{\chi}, \varphi^{\prime})$
is a pair of morphisms $A \rightarrow A^{\prime}$ (in $\CA$) and
$A_{\chi} \rightarrow A^{\prime}_{\chi}$ (in $\CA_{\chi}$) which makes
the obvious diagram commute.  There are natural forgetful functors
$\CA_{1,\chi} \rightarrow \CA$ and $\CA_{1,\chi} \rightarrow \CA_{\chi}$. 

\begin{lemma} The functors $\CA_{1,\chi} \rightarrow \CA$ and
$\CA_{1,\chi} \rightarrow \CA_{\chi}$ are equivalences of categories.
\end{lemma}
\proof
Both functors are clearly fully faithful; it suffices to
show that they are dense.  Fix an $A \in \CA$; then it is enough to
show that we have an $A^{\prime} \in \CA_{\chi}$ and a map
$\varphi: A \rightarrow A^{\prime}$ with the properties described above. 

Fix an isogeny $\phi: J \rightarrow A$.  Multiplying by a sufficiently
large integer $n$, we may assume that $\ker \phi$ contains
the kernel of the map $\psi: J \rightarrow J^{\vee}$.  Let $A^{\prime} =
J^{\vee}/\psi(\ker \phi)$.  The natural map $J \rightarrow A^{\prime}$
induced by $\psi$ has the same kernel as $\phi$, and thus descends
to an isomorphism $\varphi$ from $A$ to $A^{\prime}$.  Since $\psi$
is defined over $K$, so is $\varphi$, and we have
$\sigma \varphi = \chi(\sigma) \varphi \sigma$ since the corresponding
relation holds for $\psi$.  It follows that $(A,A^{\prime}, \varphi)$
is an object of $A_{1,\chi}$ mapping to $A$.  

The proof for the functor $\CA_{1,\chi} \rightarrow \CA_{\chi}$ is
similar.
\endproof

The upshot is that $\CA$ and $\CA_{\chi}$ are equivalent categories.
In particular the above equivalences induce a functor
$A \mapsto A_{\chi}$ from $\CA$ to $\CA_{\chi}$ and the inverse
functor $A^{\prime} \mapsto
A^{\prime}_{\chi^{-1}}$, and these functors are equivalences of
categories.  (These functors can be thought of as ``twisting'' an abelian
variety $A$ by the character $\chi$ or its inverse.  In particular
$A$ and $A_{\chi}$ become isomorphic over the splitting field of $\chi$,
via an isomorphism $\varphi$ which satisfies
$\sigma \varphi = \chi(\sigma) \varphi \sigma$ for all $\sigma$ in 
$\gal(K/\QQ)$.)  This allows us to defined a natural duality operation
on $\CA$, by taking $A^{\dual} = (A^{\vee})_{\chi^{-1}}$. 

This operation interacts with our equivalence of categories
$\CA_S \rightarrow \CM_S$ as follows:

\begin{prop} \label{prop:L}
The $\TT_S$-module $[(\Jmin)^{\dual}]^*$ is locally free of rank one.  
If $L$ denotes its inverse as a $\TT_S$-module, then we have
$[A^{\dual}] \cong [A]^* \otimes L$ for any $A \in \CA$. 
\end{prop}
\proof
The Weil pairing induces an isomorphism, $$T_{\mm} (\Jmin)^{\vee} \cong 
\Hom(T_{\mm} \Jmin, \ZZ_l(1)).$$  Since $T_{\mm} (\Jmin)^{\vee}$ and 
$T_{\mm} (\Jmin)^{\dual}$ are isomorphic as $\TT_{\mm}$-modules 
(but not Galois modules),
$T_{\mm} (\Jmin)^{\dual}$ is isomorphic to $\Hom(T_{\mm} \Jmin, \ZZ_l(1))$.
As $T_{\mm} \Jmin$ is free of rank two over $\TT_{\mm}$, it follows
that $(T_{\mm} (\Jmin)^{\dual})^*$ is free of rank two over $\TT_{\mm}$.
By Corollary~\ref{cor:tate}, $[(\Jmin)^{\dual}]^*$ is free of
rank one over $\TT_{\mm}$. 

Now for any $A \in \CA$, we have
$$
[A^{\dual}] \cong_S \Hom(A^{\dual}, \Jmin) \cong_S \Hom((\Jmin)^{\dual}, A)
\cong_S \Hom([A], [(\Jmin)^{\dual}]).$$  We have $[(\Jmin)^{\dual}] \cong
\TT_S^* \otimes L$, so $\Hom([A], [(\Jmin)^{\dual}]) \cong_S
[A]^* \otimes L$ by Lemma~\ref{lemma:star}.
\endproof
\section{The character group of the $pq$-new subvariety} \label{sec:chargp}

We now return to the situation considered in the introduction.
Fix a squarefree integer $N$, with an even number of
prime divisors, and let $\Gamma$ be a congruence subgroup
of level $M$ for some integer prime to $N$.  For $R$ a divisor of $N$,
and $D$ a divisor of $R$ which is divisible by an even number of
primes, we let $J^{D,R}$ denote the abelian variety
$J^D(\Gamma_0(N/D) \cap \Gamma)_{\RDnew}$.  Let $\TT^R$ denote
the subalgebra of $\End(J^{1,R})$ generated by the Hecke operators
$T_l$ for $l$ prime, (when $l$ divides $MN$ we take $T_l$ to be defined as
in~\cite{level}), and the
diamond bracket operators, and let $\CA^R$ denote the category of
abelian varieties defined over $\QQ$ which are $\TT^R$-equivariantly
isogenous to $J^{1,R}$.  Then for each $D$ as above, $J^{D,R}$ is
an object of $\CA^R$.  Since we will be primarily interested in
$\CA^N$, we let $J^D = J^{D,N}$ for any divisor $D$ of $N$ which is
divisible by an even number of primes.  Similarly, we take $\CA = \CA^N$,
$\TT = \TT^N$, and so forth.  Let $S^R$ be the set of maximal 
ideals of $\TT^R$ which are either Eisenstein, or of residue
characteristic 2 or 3, and take $S = S^N$.  

As mentioned in Remark~\ref{remark:modular},
the hypotheses (and hence the results) of the previous section 
hold for $\CA^R$.  We fix, once and for all, a $\Jmin_1$ inside
$\CA^1$ such that $T_{\mm} \Jmin_1$ is free of rank two over
$\TT_{\mm}^1$ for all non-Eisenstein $\mm$.  If we then take
$\Jmin_R = (\Jmin_1)_{\Rnew}$, then $T_{\mm} \Jmin_R$ is free
of rank two over $\TT^R$ for all non-Eisenstein maximal ideals $\mm$
of $\TT^R$.  We can thus consider the functors $J \rightarrow [J]_R =
\Hom(J,\Jmin_R)$ for $J \in \CA^R$, and use them to study isogenies
between the abelian varieties $J^{D,R}$ for various $D$.  By abuse of
notation, we often denote $[J]_R$ simply by $[J]$ if it is clear
from the context which category $J$ belongs to. 

$J^{1,R}$ (and hence every variety in $\CA^R$) has purely toric reduction 
at any prime $p$ dividing $R$.  Thus we can consider the character groups
of $J^{D,R}$ at such primes.  We then have relationships between
these character groups and the module $[J^{D,R}]$, given by
Proposition~\ref{prop:charisogeny}.  
Thus, if we can compute the character groups $X_p(J^{D,R})$, we can 
hope to obtain information about maps between varieties of this form from 
this data.  We carry out this program in the remainder of this paper. 

Let $R^{\prime}$ be an integer dividing $N$ and divisible by $R$.
By~\cite{CS}, Theorem 8.2,
the injection $J^{D,R^{\prime}} \rightarrow J^{D,R}$ induces an
injection on the maximal tori of the special fibers at $p$,
and hence a surjection $X_p(J^{D,R}) \rightarrow X_p(J^{D,R^{\prime}})$.  
Moreover,
both $X_p(J^{D,R})$ and $X_p(J^{D,R^{\prime}})$ are free $\ZZ$-modules, 
and there are isomorphisms: 
$$\begin{array}{ccc}
X_p(J^{D,R}) \otimes \QQ \cong \TT^R \otimes \QQ & \text{and} &
X_p(J^{D,R^{\prime}}) \otimes \QQ \cong \TT^{R^{\prime}} \otimes \QQ.
\end{array}$$  
This actually suffices
to determine the isomorphism class of $X_p(J^{D,R^{\prime}})$ in terms of 
$X_p(J^{D,R})$.

More generally, let $A$ be any finite flat $\ZZ$-algebra, and
suppose that $A \otimes \QQ$ is Gorenstein.  (This is easily
seen to hold for $A = \TT^R$ by results of the previous section.)
Let $A^{\prime}$ be a quotient of $A$ which is flat over $\ZZ$.
Let $X$ be a finitely generated $R$-module such that $X \otimes \QQ$ is
free of rank one over $A \otimes \QQ$, and $Y$ be a quotient of $X$ that
is a faithful $A^{\prime}$-module, and flat over $\ZZ$.  

\begin{lemma} \label{lemma:faithful}
Suppose $X$ and $Y$ are as above, and let $L$ be the kernel of the map
$X \rightarrow Y$.  For
an ideal $I$ of $A$, let $I^{\perp}$ be the annihilator of $I$ in $A$.
(This notation is due to Emerton~\cite{Emerton}.)
Let $K$ be the kernel of the map $A \rightarrow A^{\prime}$.
Then $L = X[K^{\perp}]$, and we have $Y \cong X \otimes_A A^{\prime}$ 
(recall our
convention, in which we take the tensor product to be the usual tensor
product modulo $\ZZ$-torsion).
\end{lemma}
\proof
We begin by showing that $K$ has finite index in $K^{\perp\perp}$.
(Note that $K^{\perp\perp}$
trivially contains $K$.)  Since $\hat A = A \otimes \QQ$ is a zero-dimensional
Gorenstein ring, we fix an isomorphism $\hat A \cong 
\Hom_{\QQ}(\hat A, \QQ)$.  For any ideal $I$ of $\hat A$, this isomorphism
identifies $I^{\perp}$ with $\Hom_{\QQ}(\hat A/I, \QQ)$.  In particular
$\dim_{\QQ} I^{\perp} = \dim_{\QQ} \hat A - \dim_{\QQ} I$, and so
$\dim_{\QQ} I^{\perp\perp} = \dim_{\QQ} I$.  Since $I$ is contained in
$I^{\perp\perp}$, the two must be equal.  It follows that $K$ is
contained in $K^{\perp\perp}$ with finite index, as required.

We have an
exact sequence:
$$0 \rightarrow L \rightarrow X \rightarrow Y \rightarrow 0.$$
Since $Y$ is an $A^{\prime}$-module, $K$ annihilates $Y$, so $KX \subset L$.
Now $KX \subset X[K^{\perp}]$, by definition of $K^{\perp}$.
Moreover, since $X \otimes \QQ$ is free over $\hat A$,
$$X[K^{\perp}] \otimes \QQ = \hat A[K^{\perp}] = K^{\perp\perp}X\otimes \QQ.$$
Thus $X[K^{\perp}]$ contains $K^{\perp\perp}X$ with finite index, and hence
contains $KX$ with finite index.  Let $n$ be the index of $KX$ in
$X[K^{\perp}]$.  Suppose $x$ in 
$X[K^{\perp}]$ maps to a nonzero element of $Y$.  Then $nx$ maps to the
zero element of $Y$, so (since we assumed $Y$ is flat over $\ZZ$), $x$ maps
to $0$ in $Y$.  Thus $X[K^{\perp}] \subset L$.

Since $X \otimes \QQ$ is generated by a single element over $A \otimes \QQ$,
$Y \otimes \QQ$ is generated by a single element over 
$A^{\prime} \otimes \QQ$; since
$Y$ is faithful over $A^{\prime}$, $Y \otimes \QQ$ is free of rank one
over $A^{\prime} \otimes \QQ$.  We thus have a commutative diagram:
$$\begin{CD}
0 @>>> L @>>> X @>>> Y @>>> 0
\\
& & @VVV @VVV @VVV \\
0 @>>> L \otimes \QQ @>>> A \otimes \QQ @>>> 
A^{\prime} \otimes \QQ @>>> 0.\\
\end{CD}$$
Since the kernel of the map $A \otimes \QQ \rightarrow A^{\prime} \otimes \QQ$
is killed by $K^{\perp}$, so is $L \otimes \QQ$.  Hence $L$ is killed by
$K^{\perp}$, so $L = X[K^{\perp}].$  In particular $L$ contains $KX$
with finite index.  Hence $Y \cong X/L$ is isomorphic to
$X/KX$, modulo $\ZZ$-torsion in $X/KX$.  The latter is precisely 
$X \otimes_A A^{\prime}$.
\endproof

We record the following useful corollary before we proceed:
\begin{cor} \label{cor:faithful} 
Let $X$, $Y$, $A$, $A^{\prime}$, and $K$ be as above. 
Then $Y^* \cong X^*[K]$.  In particular $X^*[K] \cong 
(X \otimes_A A^{\prime})^*$ for any $X$ and $K$ as above.
\end{cor}
\proof We have an exact sequence
$$0 \rightarrow L \rightarrow X \rightarrow Y \rightarrow 0.$$
Taking duals gives us the sequence:
$$0 \rightarrow Y^* \rightarrow X^* \rightarrow L^* \rightarrow 0.$$
The result now follows by the previous lemma.
\endproof

The upshot of Lemma~\ref{lemma:faithful} is the following
corollary:

\begin{cor} \label{cor:newchargp}
Let $R$ and $R^{\prime}$ be divisors of $N$, and suppose $R$ divides 
$R^{\prime}$.
For any $p$ dividing $R$, we have
$X_p(J^{D,R^{\prime}}) \cong X_p(J^{D,R}) \otimes_{\TT^R} \TT^{R^{\prime}}.$
\end{cor}

We now fix two distinct primes $p$ and $q$ dividing $N$.  Let $D$ be a divisor
of $N$, divisible by an even number of primes but {\em not} divisible by
$p$ or $q$.
The previous result, together with the following result of Ribet,
establishes
a relationship between the character groups $X_p(J^{D,Dpq})$ and
$X_q(J^{Dpq,Dpq})$: 

\begin{thm} \label{thm:ribexact}
There is a $\TT^D$-equivariant exact sequence:
$$0 \rightarrow X_q(J^{Dpq,Dpq}) \rightarrow X_p(J^{D,D}) \rightarrow 
X_p(J^D(\Gamma_0(N/Dq) \cap \Gamma))^2 \rightarrow 0.$$
\end{thm}
\proof When $D=1$, this is~\cite{level}, Theorem 4.1.  The same argument
extends to the case where $D>1$ using Buzzard's analogue of the 
Deligne-Rapoport theorem~\cite{Buz:DR}, except that Buzzard's
result assumes that the level structure is contained in $\Gamma_1(r)$
for some integer $r > 4$.  An argument that removes this restriction
can be found in the appendix.
\endproof

Dualizing, and applying Lemma~\ref{lemma:faithful}, we find that
$$X_q(J^{Dpq,Dpq})^* \cong X_p(J^{D,D})^* \otimes_{\TT^D} \TT^{Dpq}.$$
This is almost what we need; it remains to
compare $X_p(J^{D,D})^*$ with $X_p(J^{D,D})$.  We do this
via the monodromy pairing on $X_p(J^{D,D})$.  This is a
bilinear pairing: 
$$X_p(J^{D,D}) \times X_p((J^{D,D})^{\vee}) \rightarrow \ZZ,$$
for which the adjoint of $T_l$ acting on $J^{D,D}$ is $T_l^{\vee}$
acting on $(J^{D,D})^{\vee}.$

\begin{lemma} \label{lemma:rosati} 
Let $\Gamma^{\prime} = \Gamma_0(N/D) \cap \Gamma$, and
$R = NM/D$ be the level of $\Gamma^{\prime}$.
There exists a Hecke-equivariant isomorphism $\psi: J^{D,D}
\rightarrow (J^{D,D})^{\vee},$ defined over $\QQ(\mu_R)$.
Moreover, for any $\sigma \in \gal(\QQ(\mu_R)/\QQ)$,
$\sigma\psi = \chi(\sigma)\psi\sigma$, where
$\chi: \gal(\QQ(\mu_R)/\QQ) \rightarrow \TT^{\times}$ is the
character such that if $\sigma \zeta_R = \zeta_R^d$ for a 
primitive root of unity $\zeta_R$, then $\chi(\sigma) = \<d\>$.
\end{lemma}
\proof 
Let $w_R$ be the Atkin-Lehner involution on $\Gamma^{\prime} =
\Gamma_0(N/D) \cap \Gamma$; that is,
the involution of $J^{D,D} = J^D(\Gamma^{\prime})$ corresponding to the double coset
$\Gamma^{\prime}
\left(\begin{array}{cc} 0 & -t \\ R & 0\end{array} \right)\Gamma^{\prime}$.  
The discussion on page 193 of
~\cite{Sh}, together with~\cite{Sh}, Proposition 3.54, shows that
the Rosati involution induced by
the canonical polarization $\theta$ of $J^D(\Gamma^{\prime})$ 
sends an element $T_l$ of $\TT$ to $w_R T_l w_R$, for any $l$
prime to $D$.  (Strictly speaking, Shimura only covers the case
in which $D=1$ and $l$ is arbitrary, but one can just as easily
work locally at $l$ for any $l$ prime to $D$ and his argument carries
over without change.) 

For $l$ dividing $D$, one can view $X^D(\Gamma^{\prime})$ as a moduli
space of abelian surfaces with quaternionic multiplication by the
quaternion algebra $B$ of discriminant $D$, with
associated level structure. (See~\cite{Buz:DR} for details.)
From this perspective, $T_l$ is induced by the map on $X^D(\Gamma^{\prime})$
which sends an abelian surface $A$ to $A/A[\ll]$, where $\ll$
is the unique two-sided maximal ideal of index $l$ in a maximal order 
of $B_D$.  Since this map is an isomorphism, the Rosati involution
sends $T_l$ to its inverse.  But it is easy to see from the moduli
definition of $T_l$ that $T_l^2 = \<l\>$.  Thus the Rosati involution
sends $T_l$ to $\<l\>^{-1}T_l$.  Moreover, by~\cite{Sh}, Proposition
3.55, we have $w_R T_l w_R = \<l\>^{-1}T_l$ on 
$J^1(\Gamma_0(D) \cap \Gamma^{\prime})$,
and hence on $J^D(\Gamma)$ as well by Jacquet-Langlands.

Finally, it is easy to verify that both the Rosati involution and
conjugation by $w_R$ send $\<d\>$ to $\<d\>^{-1}$ for any $d$.
(For the former, simply note that $\<d\>$ is induced by an automorphism
of $X^D(\Gamma^{\prime})$; the latter is part of~\cite{Sh}, Proposition 3.55).
Since the elements $T_l$ and $\<d\>$ as $l$ and $d$ vary generate $\TT$,
we find that the Rosati involution is simply conjugation by
$w_R$.  It is then easy to check that $\psi = w_R \theta$ is Hecke-equivariant.
The second statement (and the fact that $\psi$ is defined over $\QQ(\mu_R)$),
follow from the identity $\sigma w_R = \chi(\sigma)w_R\sigma$.
\endproof

The map $\psi$ defined above induces a duality operation 
$A \rightarrow A^{\dual}$, by results of the previous section.
It also defines
an isomorphism $J^{D,D} \cong (J^{D,D})^{\vee}$
over $\QQ(\mu_R)$.
Since the character groups at $p$ of $(J^{D,D})^{\vee}$ 
are the same whether we consider the varieties 
over $\QQ$ or over $\QQ(\mu_R)$, the map $\psi$
induces an isomorphism $X_p(J^{D,D})^{\vee} \cong X_p(J^{D,D})$
(as $\TT^D$-modules but not Galois modules).
Since $\psi$ is $\TT^D$-equivariant, in the sense that $T_l^{\vee} \psi
= \psi T_l$, this identification allows us to
view the monodromy pairing 
$$X_p(J^{D,D}) \times X_p((J^{D,D})^{\vee}) \rightarrow \ZZ$$
as a $\TT^D$-equivariant bilinear
pairing on $X_p(J^{D,D})$, i.e., as a $\TT^D$-equivariant map
$X_p(J^{D,D}) \rightarrow X_p(J^{D,D})^*$.  The cokernel of
this map is naturally isomorphic to the component group of $J^{D,D}$
at $p$, and Ribet has shown: 

\begin{prop} \label{prop:component} 
The component group of $J^{D,D}$ at $p$ is supported
on $S^D$. 
\end{prop}
\proof 
Since $p$ is prime to $D$, the methods of~\cite{ribchar}, which cover 
the discriminant 1 case, carry
over to the case of arbitrary discriminant once one invokes Buzzard's analogue 
of the
Deligne-Rapoport Theorem (\cite{Buz:DR}), and works with abelian surfaces
with quaternionic multiplication rather than elliptic curves.  As
before, see the appendix for an extension of Buzzard's result
to the case in which the level structure does not contain $\Gamma_1(r)$
for some $r > 4$.\endproof

With this in hand, we can now show:

\begin{prop} \label{prop:chargp}
We have $X_p(J^{D,Dpq}) \cong_{S^{Dpq}}
X_q(J^{Dpq,Dpq})^*$.
\end{prop}
\proof
Corollary~\ref{cor:newchargp} identifies
$X_p(J^{D,Dpq})$ with $X_p(J^{D,D}) \otimes_{\TT^D} \TT^{Dpq}.$ 
On the other hand, the $\ZZ$-dual of the exact sequence of 
Theorem~\ref{thm:ribexact} identifies $X_q(J^{Dpq,Dpq})^*$ with the module
$X_p(J^{D,D})^* \otimes_{\TT^D} \TT^{Dpq}$ by Lemma~\ref{lemma:faithful}. 
Since $X_p(J^{D,D}) \cong_{S^D} X_p(J^{D,D})^*$ by
the above discussion, it follows that $X_p(J^{D,Dpq}) \cong_{S^{Dpq}} 
X_q(J^{Dpq,Dpq})^*.$
\endproof 

We will also have need of a compatibility property which the above 
isomorphisms possess.  Let $p$ and $p^{\prime}$ be two distinct primes dividing
$D$, take $\Dp = \frac{D}{pp^{\prime}}$, and choose distinct
primes $q_1$ and $q_2$ dividing $N$ but not $D$.  Let $\bq = q_1q_2.$ 
\begin{lemma} \label{lemma:compatible}
Let $f: J^{\Dp\bq,\Dp\bq} \rightarrow J^{\Dp,\Dp}$ and
$g: J^{D\bq,D\bq} \rightarrow J^{D,D}$ 
be maps such that there is a commutative diagram:
$$\begin{CD}
X_p(J^{D\bq,D\bq}) @>>> X_{p^{\prime}}(J^{\Dp\bq,\Dp\bq})\\
@Vg^{\dual}VV @VVf^{\dual}V\\
X_p(J^{D,D}) @>>> 
X_{p^{\prime}}(J^{\Dp, \Dp})
\end{CD}$$
in which the the horizontal maps
are those of Theorem~\ref{thm:ribexact}.
Let $\overline{f}$ denote the restriction of $f$ to
the $pp^{\prime}$-new subvariety $J^{\Dp \bq, D\bq}$ of
$J^{\Dp \bq, \Dp \bq}$.
Then we have a commutative diagram:
$$\begin{CD}
X_p(J^{D\bq,D\bq}) @>>> X_{p^{\prime}}(J^{\Dp\bq, D\bq})^*\\
@Vg^{\dual}VV @VV\overline{f}^*V\\
X_p(J^{D, D}) @>>> 
X_{p^{\prime}}(J^{\Dp,D})^*,
\end{CD}$$
in which the horizontal maps are the isomorphisms constructed in
Proposition~\ref{prop:chargp}. 
\end{lemma}
\proof
Consider the diagram:
$$\begin{CD}
X_p(J^{D\bq,D\bq}) @>>> X_{p^{\prime}}(J^{\Dp\bq,\Dp\bq}) @>>> 
X_{p^{\prime}}(J^{\Dp\bq,\Dp\bq})^* @<<< X_{p^{\prime}}(J^{\Dp\bq,D\bq})^*\\
@Vg^{\dual}VV @Vf^{\dual}VV @VVf^*V @VV\overline{f}^*V\\
X_p(J^{D, D}) @>>> X_{p^{\prime}}(J^{\Dp,
\Dp}) @>>> X_{p^{\prime}}(J^{\Dp, \Dp})^* 
@<<< X_{p^{\prime}}(J^{\Dp, D})^*.
\end{CD}$$
Here the left-hand square is that appearing in the hypothesis of the 
lemma.  The middle square commutes because
$f^{\dual}$ and $f$ are adjoints under the monodromy pairing.  The
right-hand square is the $\ZZ$-dual of the diagram:
$$\begin{CD}
X_{p^{\prime}}(J^{\Dp, D}) @>>> 
X_{p^{\prime}}(J^{\Dp , \Dp })\\
@VfVV @VV\overline{f}V\\
X_{p^{\prime}}(J^{\Dp\bq, D\bq}) @>>> X_{p^{\prime}}(J^{\Dp\bq,\Dp\bq}),
\end{CD}$$
which clearly commutes since $\overline{f}$ is the restriction of $f$
to the $pp^{\prime}$-new subvariety, and the horizontal maps
are just inclusions of the $pp^{\prime}$-new subvarieties into the
corresponding variety.

The isomorphisms:
$$X_p(J^{D,D}) \rightarrow X_{p^{\prime}}(J^{\Dp,D})^*,$$
$$X_p(J^{D\bq, D\bq}) \rightarrow X_{p^{\prime}}(J^{\Dp \bq,D\bq})^*$$ 
are obtained as the compositions of the two leftmost horizontal map
with the ``inverse'' of the rightmost horizontal map (bearing in mind
that the image of the rightmost horizontal map is the same as the image
of the composition of the two rightmost), in the corresponding row
of the above diagram.  The commutativity of the desired square thus
follows immediately.
\endproof

We will also need the following result about character groups in
section~\ref{sec:raise}.
\begin{prop} \label{prop:ribexact2}
Let $q_1$ and $q_2$ be two distinct primes dividing $N$ but not $D$.
Let $Y_{q_1} = X_{q_1}(J^{Dq_1q_2,Dq_1q_2})$
and let $X_{q_2} =
X_{q_2}(J^D(\Gamma_0(N/Dq_1) \cap \Gamma))^2.$
Then there is a natural $S^D$-isomorphism of $\TT^D$-modules:
$$Y_{q_1}^*/Y_{q_1} \cong_{S^D} X_{q_2}/(T_{q_1}^2 - \<q_1\>)X_{q_2}.$$
(Here $Y_{q_1}$ is considered as a submodule of $Y_{q_1}^*$ via
the monodromy pairing.)
\end{prop}
\proof In~\cite{level}, Ribet constructs an exact sequence:
$$0 \rightarrow K \rightarrow X_{q_2}/(T_{q_1}^2 - \<q_1\>)X_{q_2}
\rightarrow Y_{q_1}^*/Y_{q_1} \rightarrow 0,$$ where $K$ and $C$
are supported on $S^D$, in the case $D=1$ and $\Gamma = \Gamma_0(N)$.  
As with Theorem~\ref{thm:ribexact}, his arguments carry over to our setting 
once one invokes the analogue of the Deligne-Rapoport theorem found in
the appendix. 
\endproof

\section{Induced maps on character groups} \label{sec:maps}
In light of the character group computations of the previous section,
we might hope to study morphisms in $\CA^N$ by studying the maps they
induce on character groups.  The main tool we will use to do this is
the following result, whose proof will occupy this section as well as
the next.

\begin{thm} \label{thm:main}
Fix $D$ dividing $N$.  For any prime $p$ dividing $D$, the natural map 
$$[J^{D,D}] \rightarrow \Hom(X_p(\Jmin_D), X_p(J^{D,D}))$$ which sends
a map $\Jmin_D \rightarrow J^{D,D}$ to the induced map on character
groups at $p$ is an $S^D$-isomorphism.
\end{thm}

We will prove this theorem by constructing isomorphisms:
$$[J^{D,D}]_{\mm} \rightarrow \Hom(X_p(\Jmin_D), X_p(J^{D,D}))_{\mm}$$ 
for each $\mm$ outside $S^D$.  Although the isomorphisms we construct in
this way are in some sense natural, showing any sort of compatibility
between any two of them seems difficult.  The following lemma 
provides us a way around this difficulty:

\begin{lemma} \label{lemma:hom}
Let $M$ and $N$ be $\TT^D_{\mm}$-modules, flat over $\ZZ$, and suppose that
$M \otimes \QQ$ and $N \otimes \QQ$ are free of rank one over
$\TT^D_{\mm} \otimes \QQ$.  Suppose we have a $\TT^D_{\mm}$ module $H$,
with no $\ZZ$-torsion, and an isomorphism $f:M \otimes H \rightarrow N$.
Then the map $M \otimes \Hom(M,N) \rightarrow N$ given by evaluation
is an isomorphism.  Moreover, there is a natural map
$g:H \rightarrow \Hom(M,N)$ such that the diagram:
$$\begin{CD}
M \otimes H @>f>>  N\\
@V1 \otimes g VV @VV\Id V\\
M \otimes \Hom(M,N) @>>> N
\end{CD}$$
commutes, where the right-hand vertical map is the identity.  Finally,
if $H \cong \Hom(M,N)$, then $g$ is an isomorphism.
\end{lemma}
\proof
The natural map $g$ takes an element $h$ of $H$ to  
the map $m \mapsto f(h \otimes m).$  The commutativity of the above
diagram is trivial.  Moreover, the upper horizontal map is an isomorphism
by hypothesis, and the right-hand vertical map is the identity.  It follows 
that the lower horizontal map
is surjective, and hence (by counting ranks over $\ZZ$) an isomorphism.
(Recall that $M \otimes \Hom(M,N)$ is torsion free by convention.)

It remains to prove that if $H \cong \Hom(M,N)$, then the natural map
given above is an isomorphism.  Composing the inverse of the isomorphism
$H \rightarrow \Hom(M,N)$ with $g$,
we obtain an endomorphism $\sigma$ of $\Hom(M,N)$, such that 
$g(H) \subset \Hom(M,N)$ is equal to
$\sigma\Hom(M,N)$.  It thus suffices to show that $\sigma$ is a unit
in $\End(\Hom(M,N))$.  The map $M \otimes \Hom(M,N) \cong N$ induces
an inclusion of $\End(\Hom(M,N))$ into $\End(N)$.  It follows that the
image of $M \otimes H$ in $N$, under the composition of the left-hand
vertical map and the bottom horizontal map, is equal to $\sigma N$.
By the commutativity of the diagram, $\sigma$ is a unit in $\End(N)$.
Suppose $\sigma$ were not a unit in $\End(\Hom(M,N))$.  Then there
would be a maximal ideal $\mm_1$ of $\End(\Hom(M,N))$ containing
$\sigma$.  Since $\End(\Hom(M,N))$ is contained in $\End(N)$ with finite
index, there is a maximal ideal $\hat \mm_1$ of $\End(N)$ containing $\mm_1$.
In particular, $\hat \mm_1$ would contain $\sigma$, 
and so $\sigma$ would not be a unit in $\End(N)$.  This is a
contradiction, so $\sigma$ is a unit in $\End(\Hom(M,N))$, as required.
\endproof

The usefulness of this lemma can be seen in the following
corollary:
\begin{cor} \label{cor:hom}
Let $J \in \CA^D$, $p$ be a prime dividing $R$, $\mm$ a maximal ideal
of $\TT^R$ outside $S_D$, and
suppose that $[J]_{\mm} \cong \Hom(X_p(\Jmin_D),X_p(J))_{\mm}$.
Then the natural map $[J]_{\mm} \rightarrow \Hom(X_p(\Jmin_D),X_p(J))_{\mm}$
is an isomorphism.
\end{cor}
\proof
We have a natural isomorphism: $X_p(\Jmin_D)_{\mm} \otimes [J]_{\mm}
\rightarrow X_p(J)_{\mm}.$  By Lemma~\ref{lemma:hom}, this
induces a map $[J]_{\mm} \rightarrow \Hom(X_p(\Jmin_D),X_p(J))_{\mm}$,
which is easily seen to be the natural map $[J]_{\mm} \rightarrow
\Hom(X_p(\Jmin_D), X_p(J))_{\mm}.$  Since the two modules are isomorphic,
this map is an isomorphism by Lemma~\ref{lemma:hom}.
\endproof

In particular, in order to prove Theorem~\ref{thm:main}, it suffices to
construct, for each $\mm$ outside $S^D$, an arbitrary isomorphism
$$[J^{D,D}]_{\mm} \rightarrow \Hom(X_p(\Jmin_D),X_p(J^{D,D}))_{\mm}.$$
We begin by doing this for maximal ideals $\mm$ with a certain property.

\begin{defn} Suppose $D$ divides $N$, and is divisible by an even number
of primes, and 
let $\mm$ be a maximal ideal of $\TT^D$ such that the
representation $\overline{\rho}_{\mm}$ is absolutely irreducible.
Then $\mm$ is {\em controllable} at $p$ (for some $p$ dividing $N$)
if one of the following conditions holds:
\begin{enumerate}
\item $\overline{\rho}_{\mm}$ is not finite at $p$, or
\item $\overline{\rho}_{\mm}$ is unramified at $p$, $p \neq l$, and 
$\overline{\rho}_{\mm}(\Frob_p)$ is not a scalar, or
\item $p = l$, and $l \ne 2$, or
\item $p = l = 2$, and the restriction of 
$\overline{\rho}_{\mm}$ to a decomposition group at $2$ is not
contained in the scalar matrices.
\end{enumerate}
\end{defn}

The point of this definition is the following lemma:

\begin{lemma} \label{lemma:control} 
Suppose that $\mm$ is controllable at $p$. Then $X_p(\Jmin_D)_{\mm}$ 
is free of rank one over $\TT^D_{\mm}$.
\end{lemma}

\proof 
This is essentially Mazur's principle.
By Nakayama's lemma, and the fact that $X_p(\Jmin_D)_{\mm}$ is faithful
over $\TT^D_{\mm}$, it suffices to show that the quotient 
$X_p(\Jmin_D)/\mm X_p(\Jmin_D)$
has dimension one.  We will prove the lemma in this form.

Let $\Jmin_D[\mm]^f$ denote the maximal $\TT^D$-stable subgroup of
$\Jmin_D[\mm]$ which is finite at $p$- i.e., which extends to a finite, flat
group scheme $W$ over $\ZZ_p$.  Then we have an injection
$W_{\FF_p} \rightarrow (\Jmin_D)_{\FF_p}[\mm]$, where $(\Jmin_D)_{\FF_p}$
is the special fiber of the Ner{\'o}n model of $\Jmin_D$ over $\FF_p$. 
Since $W$ is flat over $\ZZ_p$, the dimension of $W_{\FF_p}$ and
the dimension of the general fiber of $W$ over $\TT/\mm$ agree.
Since the latter is at most two, and is exactly two if and only if
$\Jmin_D[\mm]$ is finite at $p$, it follows that
$(\Jmin_D)_{\FF_p}[\mm]$ has dimension at most two, with equality only if
$(\Jmin_D)[\mm]$ is finite at $p$.

Now, since $(\Jmin_D)_{\FF_p}$ is an extension of the torus
$\Hom(X_p(\Jmin_D), \GG_m)$ by a finite group, and  
$$\dim X_p(\Jmin_D)/\mm X_p(\Jmin_D) \geq 1,$$ we have
$$1 \leq \dim X_p(\Jmin_D)/ \mm X_p(\Jmin_D) \leq 2,$$ with equality on
the right if and only if $\Jmin_D[\mm]$ is finite at $p$.  In particular,
if $\overline{\rho}_{\mm} = \Hom(\Jmin_D[\mm], \mu_l)$ is not finite
at $p$, then $\Jmin_D[\mm]$ is not finite at $p$, and so
$X_p(\Jmin_D)/ \mm X_p(\Jmin_D)$ has dimension one.

We may thus assume that $\overline{\rho}_{\mm}$, and hence $\Jmin_D[\mm]$,
is finite at $p$.  Suppose further that $X_p(\Jmin_D)/\mm X_p(\Jmin_D)$
has dimension two.  Then we are in the following situation:
$\Jmin_D[\mm]$ extends to a finite flat group scheme $W$ over $\ZZ_p$,
and the special fiber of $W$ over $\FF_p$ is the multiplicative group scheme 
$\Hom(X_p(\Jmin_D)/\mm X_p(\Jmin_D), \mu_l).$  In particular,
the Cartier dual $W^{\vee}$ of $W$ is an \'etale group scheme over
$\ZZ_p$.  Since $\gal(\overline{\QQ}/\QQ)$ acts on the general fiber
of $W^{\vee}$ by $\overline{\rho}_{\mm}$, it follows that 
$\overline{\rho}_{\mm}$ is
unramified at $p$.  In particular, the image of a decomposition group
at $p$ under $\overline{\rho}_{\mm}$ is determined by the action of $\Frob_p$
on the special fiber of $W^{\vee}$ = $X_p/\mm X_p$.
 
Suppose first that $p \neq l$.  Then $\Frob_p$ acts on 
$X_p(\Jmin_D)/\mm X_p(\Jmin_D)$ via the Atkin-Lehner operator 
$w_p^{-1}$, as $X_p(\Jmin_D)$
can be identified with a subspace of $X_p(J^{1,D})$, on which this
is well-known.
Since $w_p$ acts as the scalar
$-T_p$ on this space, it follows that $\overline{\rho}_{\mm}(\Frob_p)$
is a scalar.  In particular, $\mm$ cannot be controllable at $p$ in this
case.

Suppose finally that $p=l$.  
Since the determinant of $\overline{\rho}_{\mm}$ is the mod $l$
cyclotomic character times an even character of finite order, 
$\overline{\rho}_{\mm}$ is ramified at $p$ unless $p=l=2$.  This
contradicts the fact that $W^{\vee}$ is {\'e}tale.  If $p=l=2$,
then we find as before that the action of $\Frob_p$ on
$W^{\vee} = X_p/\mm X_p$ is by scalars, and hence $\mm$ again
cannot be controllable at $p$. 
\endproof

When this criterion is satisfied for some prime dividing $D$, we
may easily construct the desired isomorphism:

\begin{cor} \label{cor:key} 
Let $p$ divide $D$, and suppose $\mm$ is
not Eisenstein and controllable at at least one prime dividing $D$.  
Then $[J^{D,D}]_{\mm}$ is isomorphic to  
$\Hom(X_p(\Jmin_D),X_p(J^{D,D}))_{\mm}$ for any $\mm$ outside $S^D$.
\end{cor}
\proof
If $\mm$ is controllable at $p$, then $X_p(\Jmin_D)_{\mm}$ is free of rank
one over $\TT_{\mm}$, and $X_p(J^{D,D})_{\mm} \cong [J^{D,D}]_{\mm} \otimes 
X_p(\Jmin_D)_{\mm}.$  The result follows immediately.

If $\mm$ is not controllable at $p$, there is another prime $q$ dividing $D$
at which $\mm$ is controllable.  Then
we have a sequence of isomorphisms:
$$\begin{array}{rcll}
\Hom(X_p(\Jmin_D),X_p(J^{D,D}))_{\mm} & \cong &
\Hom(X_p(\Jmin_D),X_q(J^{\frac{D}{pq},D})^*)_{\mm} &
\mbox{(Proposition~\ref{prop:chargp})}\\
& \cong & (X_p(\Jmin_D) \otimes X_q(J^{\frac{D}{pq},D}))^*_{\mm} &
\mbox{(Lemma~\ref{lemma:star})} \\
& \cong & (X_p(\Jmin_D) \otimes [J^{\frac{D}{pq},D}])^*_{\mm} &
\mbox{(controllability at $q$)} \\
& \cong & X_p(J^{\frac{D}{pq},D})_{\mm}^* &
\mbox{(Proposition~\ref{prop:charisogeny})} \\
& \cong & X_q(J^{D,D})_{\mm} &
\mbox{(Proposition~\ref{prop:chargp})} \\
& \cong & [J^{D,D}]_{\mm} &
\mbox{(controllability at $q$).}
\end{array}$$
This proves the first statement; the second statement follows by
Corollary~\ref{cor:hom}.
\endproof

It remains to deal with those maximal ideals $\mm$ that are not
controllable at {\em any} prime dividing $D$.  We do so in the next section,
by ``raising the level''.  Specifically, we introduce two auxiliary primes
$q_1$ and $q_2$, prime to $NM$, such that $\mm$ lifts to a maximal ideal
of $\tilde \mm$ of $\TT^D(\Gamma_0(q_1q_2N/D) \cap \Gamma)$
which is {\em new} at $q_1$ and $q_2$ and {\em controllable} at these
primes.  Then $\tilde \mm$ descends to a maximal ideal of
$\TT^{Dq_1q_2}(\Gamma_0(N/D) \cap \Gamma)$ (because $\tilde \mm$ is new)  
and an analogue of Corollary~\ref{cor:key} holds for the abelian variety
$J^{Dq_1q_2}(\Gamma_0(N/D) \cap \Gamma)$ at $\tilde \mm$, as
$\tilde \mm$ is controllable at $q_1$.  We then exploit a geometric
relationship between this abelian variety and $J^{D,D}$ to prove the
desired result for $J^{D,D}$.

\section{Raising the level} \label{sec:raise}

The key to the level raising argument is the following result:

\begin{lemma} \label{lemma:raise} 
Let $D$ be the product of an even number of distinct primes,
and $\Gamma$ a congruence subgroup of level $M$ prime to $D$.  Let
$\mm$ be a non-Eisenstein maximal ideal of the Hecke algebra $\TT^D(\Gamma)$,
of residue characteristic $l > 3$.  
Then there exists a prime $q$, not dividing $2NDl$,
such that:
\begin{enumerate}
\item $\overline{\rho}_{\mm}(\Frob_q)$ is not a scalar, and
\item There is a maximal ideal $\tilde \mm$
of $\TT^D(\Gamma_0(q) \cap \Gamma)^{\qnew}$ such that the representations
$\overline{\rho}_{\mm}$ and $\overline{\rho}_{\tilde \mm}$ are isomorphic.
\end{enumerate}
\end{lemma}
\proof
Let $\chi$ be the character such that $\det \overline{\rho}_{\mm} = \chi \omega$,
where $\omega$ is the cyclotomic character.  Let $\sigma$ be the
representation $\overline{\rho}_{\mm} \times \chi$.
By the {\v C}ebotarev density theorem, we may choose $q$ such that
$\sigma(\Frob_q)$ is conjugate to $\sigma(c)$, where $c$ denotes complex 
conjugation.  Arguments of Diamond~\cite{Dlevel} show that
$\overline{\rho}_{\mm}(\Frob_q)$ then has order two, trace zero and determinant 
$-1$.  
In this setting (1) is immediate; (2) follows from
~\cite{DT}, Theorem A.
\endproof

Suppose we fix a $q$ and $\tilde \mm$ as above.  Then $\tilde \mm$
is $q$-old; that is, it descends to an ideal of
the algebra $\TT^D(\Gamma_0(q) \cap \Gamma)^{\qold}.$  We
now investigate the structure of this ring.  The two degeneracy maps 
$$J^D(\Gamma) \rightarrow J^D(\Gamma_0(q) \cap \Gamma)$$ induce
an isogeny $$J^D(\Gamma)^2 \rightarrow J^D(\Gamma_0(q) \cap \Gamma)_{\qold}.$$ 
This isogeny is compatible with a certain action of the Hecke operators
on $J^D(\Gamma)^2$, in which the diamond bracket operators and the operators
$T_n$ for $n$ prime to $q$ act diagonally, and the action of $T_q$
is given by the matrix 
$$\left( \begin{array}{cc}\tau_q & q\<q\> \\ -1 & 0 \end{array}\right).$$  
(Here $\tau_q$ denotes the endomorphism of $J^D(\Gamma)^2$ induced by the
$q$th Hecke operator acting diagonally.)  A proof may
be found in~\cite{Dlevel}.
 
It follows that $\TT^D(\Gamma_0(q) \cap \Gamma)^{\qold}$ is isomorphic
to $R[T_q]/(T_q^2 - \tau_q T_q + q\<q\>)$, where $R$ is the subalgebra
of $\TT^D(\Gamma)$ generated by the diamond bracket operators and the
Hecke operators $T_n$ for $n$ prime to $q$.
By the lemma on p. 491 of~\cite{wiles}, $R$ is all of 
$\TT^D(\Gamma)$.  Moreover, the Eichler-Shimura relation implies that
$\tilde \mm \cap R = \mm.$

\begin{prop} \label{prop:fusion}
Let $q$ and $\tilde \mm$ be chosen as above.
Then 
$$\dim J^D(\Gamma_0(q) \cap \Gamma)[\tilde \mm] =
\dim J^D(\Gamma_0(q) \cap \Gamma)_{\qnew}[\tilde \mm] = 
\dim J^D(\Gamma_0(q) \cap \Gamma)_{\qold}[\tilde \mm].$$
\end{prop}
\proof
To make the notation
more concise, let $\tilde J = J^D(\Gamma_0(q) \cap \Gamma)$,
$\tilde J_{\new} = \tilde J_{\qnew}$, and $\tilde J_{\old} = \tilde J_{\qold}$.
Let $J = J^D(\Gamma)$.  As described above, the two degeneracy maps 
$J \rightarrow \tilde J$
induce a Hecke-equivariant isogeny $(J^2) \rightarrow \tilde J_{\old}$.  
By~\cite{DT}, Theorem 2,
the kernel of this isogeny is Eisenstein; in particular it is supported
away from $\tilde \mm$. 

Let $\tilde J^{\old}$ be the $q$-old {\em quotient} of $\tilde J$
(and thus the dual of $\tilde J_{\old}$), and
consider the sequence of maps:
$$J^2 \rightarrow \tilde J_{\old} \stackrel{i}{\rightarrow} \tilde J 
\stackrel{\psi}{\rightarrow}
(\tilde J)^{\dual} \stackrel{i^{\dual}}{\rightarrow} 
\tilde J^{\old} \rightarrow J^2,$$ where the map $\psi$
is given by Lemma~\ref{lemma:rosati}, and the map $J^{\old} \rightarrow J^2$
is the dual of the isogeny $J^2 \rightarrow \tilde J_{\old}$.  Calculations
of Ribet~\cite{level} 
show that the composition of these maps is
given by the matrix $\left(\begin{array}{cc}q+1 & \<q\>T_q \\ T_q & q+1
\end{array} \right)$, and is
in particular equal to $T_q^2 - \<q\>$ times an automorphism of $J^2$.
Let $\psi_{\old}: \tilde J_{\old} \rightarrow \tilde J^{\old}$ be
the map $i^{\dual}\psi i$, above.  The kernel of $\psi_{\old}$ is
equal to the intersection of $\tilde J_{\old}$ with the kernel of the
map $\tilde J \rightarrow \tilde J^{\old}$; the latter is precisely
$\tilde J_{\new}$.  Since the intersection of $\tilde J_{\old}$ with
$\tilde J_{\new}$ is finite, and $\tilde J^{\old}$ and $\tilde J_{\old}$ 
have the same dimension, $\psi_{\old}$ is surjective.
Since the map $J^2 \rightarrow \tJ_{\old}$ is an isogeny,
and $T_q^2 - \<q\>: J^2 \rightarrow J^2$ factors as 
$$J^2 \rightarrow \tJ_{\old} \stackrel{\psi_{\old}}{\rightarrow}
\tJ^{\old} \rightarrow J^2,$$
it follows that the map $T_q^2 - \<q\>$ is a surjective endomorphism
of $J^2$.  We immediately see that the map 
$T_q^2 - \<q\>: \tilde J_{\old} \rightarrow
\tilde J_{\old}$ is surjective as well.  Finally, observe that
since the kernel of the map $J^2 \rightarrow \tilde J_{\old}$ is supported
away from $\tilde \mm$, it follows that the $\tilde \mm$-primary part of 
$\ker \psi_{\old}$ is equal to the 
$\tilde \mm$-primary part of $J_{\old}[T_q^2 - \<q\>]$.  

The operator $T_q$ acts on $\tilde J_{\new}$ as $-w_q$, where 
$w_q$ is the Atkin-Lehner operator.  Since $w_q^2 = \<q\>$, it follows
that $T_q^2 - \<q\>$ annihilates $\tilde J_{\new}$.  The endomorphism
of $\tilde J$ induced by this map thus factors through 
$\tilde J^{\old}$; its image therefore lies in
the $q$-old subvariety $\tilde J_{\old}$.  Since 
$(T_q^2-\<q\>)\tilde J_{\old}$ is equal to $\tilde J_{\old}$, this
image must be all of $\tilde J_{\old}$. 

Now consider the diagram:
$$\begin{CD}
0 @>>> \tilde J_{\old} \cap \tilde J_{\new} @>>> 
\tilde J_{\old} \oplus \tilde J_{\new} @>>> \tilde J 
@>>> 0\\
& & @VT_q^2-\<q\>VV @VT_q^2-\<q\>VV @VT_q^2-\<q\>VV\\
0 @>>> \tilde J_{\old} \cap \tilde J_{\new} @>>> 
\tilde J_{\old} \oplus \tilde J_{\new} @>>> \tilde J 
@>>> 0.
\end{CD}$$
The Snake lemma then yields an exact sequence:
$$0 \rightarrow \tilde J_{\old} \cap \tilde J_{\new} \rightarrow
\tilde J_{\old}[T_q^2 - \<q\>] \oplus \tilde J_{\new} \rightarrow
\tilde J[T_q^2 - \<q\>] \rightarrow \tilde J_{\old} \cap \tilde J_{\new}
\rightarrow \tilde J_{\new}.$$
The final map in this sequence is simply the inclusion
$\tilde J_{\old} \cap \tilde J_{\new} \rightarrow \tilde J_{\new}$,
and hence injective, so this yields a short exact sequence:
$$0 \rightarrow \tilde J_{\old} \cap \tilde J_{\new} \rightarrow
\tilde J_{\old}[T_q^2 - \<q\>] \oplus \tilde J_{\new} \rightarrow
\tilde J[T_q^2 - \<q\>] \rightarrow 0.$$  Taking $\tilde \mm^{\infty}$-torsion
then yields
$$0 \rightarrow (\tilde J_{\old} \cap \tilde J_{\new})_{\tilde \mm}
\rightarrow \tilde J_{\old}[T_q^2 - \<q\>]_{\tilde \mm} \oplus 
\tilde J_{\new}[\tilde \mm^{\infty}] \rightarrow \tilde 
J[T_q^2 - \<q\>, \tilde \mm^{\infty}]
\rightarrow 0.$$ 

The discussion above showed that $J_{\old}[T_q^2 - \<q\>]_{\tmm}$
lies in the kernel of $\psi_{\old}$, and hence lies in
$J_{\new}[T_q^2 - \<q\>]$.  It follows that
$(J_{\old} \cap J_{\new})_{\tmm} = J_{\old}[T_q^2 - \<q\>]_{\tmm}$,
and hence that $\tilde J_{\new}[\tilde \mm^{\infty}]$ is isomorphic to
$\tilde J[T_q^2-\<q\>, \tilde \mm^{\infty}]$.  Since $T_q^2 - \<q\>$
lies in $\tilde \mm$ (since $\tilde \mm$ is $q$-new), it follows
that $\tilde J_{\new}[\tilde \mm] = \tilde J[\tilde \mm]$.

It remains to show that $\tilde J_{\old}[\tilde \mm] = \tilde J[\tilde \mm].$ 
Suppose $\tilde J[\tilde \mm]$ has dimension $2n.$
Consider the connected component $\tilde J_q^0$
of the special fiber of the Ner\'on model of $\tilde J$ at $q$.  Since
$\tilde J$ has semistable reduction at $q$, it fits into an exact sequence
$$0 \rightarrow \Hom(X_q(\tilde J), \GG_m) \rightarrow \tilde J_q^0
\rightarrow \tilde J^2 \rightarrow 0.$$  Since the natural isogeny
$J^2 \rightarrow \tilde J_{\old}$ has kernel supported away from $\mm$,
it suffices to show that $J^2[\tilde \mm]$ has dimension $2n$.
 
The new subvariety $\tilde J_{\new}$ has purely toric reduction at $q$,
so the special fiber of its Ner\'on model at $q$ is an extension of
the torus $\Hom(X_q(\tilde J_{\new}), \GG_m)$ by the finite group
$\Phi_q(\tilde J_{\new})$.  
As $\overline{\rho}_{\mm}$ is unramified at $q$, and $\tilde J_{\new}[\tilde \mm]$
has dimension $2n$, $(\tilde J_{\new})_q[\tilde \mm]$
has dimension $2n$.  We have an exact sequence:
$$0 \rightarrow \Hom(X_q(\tilde J_{\new}), \GG_m)[\tilde \mm]
\rightarrow (\tilde J_{\new})_q[\tilde \mm] \rightarrow 
\Phi_q(\tilde J_{\new})[\tilde \mm].$$ 
Moreover, since $\tilde \mm$ is controllable at $q$ ($\overline{\rho}_{\mm}(\Frob_q)$
is not a scalar), it follows by Lemma~\ref{lemma:control} and
Corollary~\ref{cor:controllable} that $T_{\tilde \mm} \tilde J_{\new} \cong
X_q(\tilde J_{\new})_{\tilde \mm}^2.$  In particular, $X_q(\tilde J_{\new})/
\tilde \mm X_q(\tilde J_{\new})$ has dimension $n$.  The above exact
sequence then shows that $\Phi_q(\tilde J_{\new})[\tilde \mm]$ has
dimension at least $n$.  (We will see later that it must have dimension
exactly $n$.)

The inclusion $\tilde J_{\new} \rightarrow
\tilde J$ induces an injection on the special fibers $(J_{\new})_q$
and $\tilde J_q$ of the corresponding Ner\'on models,
and in particular an isomorphism $X_q(\tilde J_{\new}) \rightarrow X_q(J)$.
Since the group of connected
components $\Phi_q(\tilde J)$ is Eisenstein (c.f. 
Proposition~\ref{prop:component}),
$\Phi_q(\tilde J_{\new})[\tilde \mm]$ maps into the connected component of
$\tilde J_q$.  Since the map on special fibers is injective,
the image of $\Phi_q(\tilde J_{\new})[\tilde \mm]$ in $\tilde J_q$
has trivial intersection with the torus, and hence maps injectively
into $J^2$.

We thus obtain an injection of $\Phi_q(\tilde J_{\new})[\tilde \mm]$ 
into $J^2[\tilde \mm]$.  The latter is a direct sum of two dimensional
subspaces $V_i$ on which Frobenius acts via $\overline{\rho}_{\mm}(\Frob_q)$.
Since $\Frob_q$ acts as a scalar on $\Phi_q(\tilde J_{\new})$,
the image of $\Phi_q(\tilde J_{\new})$ in each $V_i$ has dimension at most
one.  As $\Phi_q(\tilde J_{\new})[\tilde \mm]$ has dimension at least
$n$, it follows that $J^2[\tilde \mm]$ has dimension at least $2n$.
But $\dim J^2[\tilde \mm] = \dim \tilde J_{\old}[\tilde \mm]$,
and the latter is less than $\dim \tilde J[\tilde \mm]$ which is $2n$. 
Thus $\dim \tilde J_{\old}[\tilde \mm] = \dim \tilde J[\tilde \mm]$, as
required.  \endproof

We now return to the setting of the previous section;
we have fixed a $D$ dividing $N$ which is the product of
an even number of primes, a maximal ideal $\mm$, and a prime $p$
dividing $D$, and we wish to show that 
$[J^{D,D}]_{\mm} \cong \Hom(X_p(\Jmin),X_p(J^{D,D}))_{\mm}.$
Let $p^{\prime}$ be a prime dividing $D$ other than $p$,
and let $D^{\prime} = \frac{D}{pp^{\prime}}$.

Now we choose two primes $q_1$ and $q_2$ as in Lemma~\ref{lemma:raise}.
The choice of these primes determines several abelian varieties.  In
particular, the abelian varieties $J^{\Dp}(\Gamma_0(q_1q_2N/\Dp) \cap
\Gamma)$ and  
$J^D(\Gamma_0(q_1q_2N/D) \cap \Gamma)$
arise from $J^{\Dp,\Dp}$ and $J^{D,D}$ by raising the
level at $q_1$ and $q_2$.  We denote them by $\tJ^{\Dp,\Dp}$
and $\tJ^{D,D}$, respectively.  We will also need to consider the
varieties $\tJ^{\Dp q_1q_2, \Dp q_1q_2} = 
J^{\Dp q_1q_2}(\Gamma_0(N/\Dp) \cap \Gamma)$ and $\tJ^{Dq_1q_2,Dq_1q_2}
= J^{Dq_1q_2}(\Gamma_0(N/D) \cap \Gamma).$  Let $\tA^{\Dp}, \tA^D,
\tA^{\Dp q_1q_2},$ and $\tA^{Dq_1q_2}$ denote the categories of
abelian varieties over $\QQ$ which are Hecke-equivariantly isogenous
to $\tJ^{\Dp,\Dp}, \tJ^{D,D}, \tJ^{\Dp q_1q_2, \Dp q_1q_2}$, and 
$\tJ^{Dq_1q_2, Dq_1q_2}$, respectively.
Also consider the corresponding categories
$(\tA^{\Dp})^{\old}$ and $(\tA^D)^{\old}$ associated to the
$q_1q_2$-old quotients of $\tJ^{\Dp}$ and $\tJ^D$, respectively.
We adopt all of the notational conventions introduced in
section~\ref{sec:chargp} for these varieties as well; for instance,
$\tJ^{D, Dq_1q_2}$ will denote the $q_1q_2$-new subvariety of
$\tJ^{D, D}$.  Finally, let $\tmm$ be a maximal ideal of
$\tTT^{\Dp}$ whose associated Galois representation is isomorphic to 
$\overline{\rho}_{\mm}$.

The four degeneracy maps $\phi_1, \phi_{q_1}, \phi_{q_2},$ and $\phi_{q_1q_2}$
from $J^{\Dp,\Dp}$ to $\tJ^{\Dp,\Dp}$ induce an isogeny 
$(J^{\Dp,\Dp})^4 \rightarrow (\tJ^{\Dp, \Dp})_{\old}$, where the subscript
$\old$ denotes the $q_1q_2$-old subvariety.  The action of $T_{q_1}$ and 
$T_{q_2}$ on $(J^{\Dp,\Dp})^4$ is given by the matrices:
$$\begin{array}{ccc}T_{q_1} = \left( 
\begin{array}{cccc}
\tau_{q_1} & q_1\<q_1\> & 0 &0\\
-1 & 0 & 0 & 0\\
0 & 0 & \tau_{q_1} & q_1\<q_1\>\\
0 & 0 & -1 & 0\\
\end{array} \right); & \quad & T_{q_2} =
\left(
\begin{array}{cccc}
\tau_{q_2} & 0 & q_2\<q_2\> & 0\\
0 & \tau_{q_2} & 0 & q_2\<q_2\>\\
-1 & 0 & 0 & 0\\
0 & -1 & 0 & 0\\
\end{array} \right),\end{array}$$
where $\tau_{q_1}$ and $\tau_{q_2}$ denote the endomorphisms
of $J^{\Dp,\Dp}$ induced by the $q_1$th and $q_2$th
Hecke operators.  The diamond bracket operators, and the Hecke operators $T_n$
for $n$ prime to $q_1q_2$, act diagonally.  
By the same arguments as those that follow the proof of 
Lemma~\ref{lemma:raise}, we have 
$$(\tTT^{\Dp})^{\old} \cong
\TT^{\Dp}[T_{q_1},T_{q_2}]/(T_{q_1}^2 - \tau_{q_1}T_{q_1} + q_1\<q_1\>,
T_{q_2}^2 - \tau_{q_2}T_{q_2} + q_2\<q_2\>).$$  In particular
$(\tTT^{\Dp})^{\old}$ is a free $\TT^{\Dp}$-algebra of rank 4. 

More generally, given any abelian variety $J$ in $\CA^{\Dp}$
or $\CA^D$, we can form the abelian variety $J^4$ on which
$(\TT^{\Dp})^{\old}$ acts via the matrices given above.  In this situation,
we have:
\begin{lemma} Let $A_1$ and $A_2$ be abelian varieties in
$\CA^{\Dp}$ or $\CA^D$.  Then 
$$\Hom(A_1^4,A_2^4) \cong 
\Hom(A_1,A_2) \otimes_{\TT^{\Dp}} (\tTT^{\Dp})^{\old},$$ 
where $\Hom$ denotes 
morphisms in $(tA^{\Dp})^{\old}$ or $(\tA^D)^{\old}$, respectively.
\end{lemma}
\proof
A map $A_1^4 \rightarrow A_2^4$ is given by a four by four matrix
of maps $A_1 \rightarrow A_2$.  For such a map to be $(\tTT^{\Dp})^{\old}$-
equivariant, it must commute with the matrices defining the action of
$T_{q_1}$ and $T_{q_2}$.  It is then straightforward to check that
every such matrix can be written in the form
$a^4 + T_{q_1}b^4 + T_{q_2}c^4 + T_{q_1}T_{q_2}d^4$ for maps $a,b,c,$ and $d$
from $A_1$ to $A_2$.  (Here $a^4: (A_1)^4 \rightarrow (A_2)^4$ is the map
obtained by evaluating the map $a$ at each entry of the four-tuple
defining a point of $(A_1)^4$.)  
That all maps of the above form commute with $T_{q_1}$ and $T_{q_2}$
is clear.
\endproof

Recall that we have fixed abelian varieties $\Jmin_{\Dp}$ and
$\Jmin_D$ in $\CA^{\Dp}$ and $\CA^D$ which satisfy ``multiplicity one''
at all maximal ideals outside $S^{\Dp}$ or $S^D$.  We will use the above
results to make ``compatible'' choices of varieties $\Jmin$ for
each of the other categories involved.  The key is the following lemma:
\begin{lemma} Let $\mm$ be a maximal ideal of $\TT^{\Dp}$ or
$\TT^D$, and $\tmm$ a maximal ideal of $(\tTT^{\Dp})^{\old}$
(or $(\tTT^D)^{\old}$) above $\mm$.  Then
$\dim_{\TT^{\Dp}/\mm} J[\mm] = \dim_{(\tTT^{\Dp})^{\old}/\tmm} (J^4)[\tmm]$
(and analogously for $\TT^D$ and $(\tTT^D)^{\old}$.)
\end{lemma} 
\proof 
The set $\{1, T_{q_1}, T_{q_2}, T_{q_1}T_{q_2}\}$ is a basis for
$(\tTT^\Dp)^{\old}$ over $\TT^\Dp$.  This basis gives rise to a dual basis
for $\Hom_{\TT^\Dp}((\tTT^\Dp)^{\old}, \TT^\Dp)$, and the matrices giving
the action of $T_{q_1}$ and $T_{q_2}$ with respect to this dual basis
are the same as the matrices giving the action of $T_{q_1}$ and
$T_{q_2}$ on $(J^\Dp)^4$.  It follows that $(J^\Dp)^4[\mm]$ is isomorphic
to $\Hom_{\TT^\Dp}((\tTT^\Dp)^{\old}, J^\Dp[\mm])$.  Therefore,
$(J^\Dp)^4[\tmm]$ is isomorphic to
$\Hom_{\TT^\Dp/\mm}((\tTT^\Dp)^{\old}/\tmm, J^\Dp[\mm])$.  The result follows
immediately, and the proof for $J^D$ is exactly the same.
\endproof

It follows from this that $(\Jmin_{\Dp})^4$ satisfies ``multiplicity one''
at every non-Eisenstein maximal ideal of $(\tTT^{\Dp})^{\old}$ of
residue characteristic greater than 3.  By Theorem~\ref{thm:equiv}, it
is straightforward to check that there is an abelian variety
$\tJmin_{\Dp}$ of $\tA^{\Dp}$ such that $\tJmin_{\Dp}$ satisfies
``multiplicity one'' at every non-Eisenstein maximal ideal of $\tTT^{\Dp}$
of residue characteristic greater than 3, and such that $(\tJmin_{\Dp})_{\old} =
(\Jmin_{\Dp})^4$.  As usual, we take $\tJmin_D = 
(\tJmin_{\Dp})_{\ppnew},$ and similarly for $\tJmin_{\Dp q_1q_2}$ and
$\tJmin_{Dq_1q_2}$.  Finally, to $(\tA^{\Dp})^{\old}$ and
$(\tA^D)^{\old}$ we associate $\tJmin_{\Dp,\old} = (\Jmin_{\Dp})^4$ and
$\tJmin_{D,\old} = (\Jmin_D)^4$, respectively.  For each category 
$\tA$ in the set
$$\{\tA^{\Dp}, \tA^D, \tA^{\Dp q_1q_2}, \tA^{Dq_1q_2}, (\tA^{\Dp})^{\old},
(\tA^D)^{\old}\}$$ the above choice
of abelian variety gives a functor $[-]_{\tA}$ from $\tA$ to
the category of ``rank one'' modules over the corresponding Hecke algebra. 

\begin{lemma} In this setting we have:
\begin{enumerate}
\item Let $\tA = \tA^{\Dp}$ or $\tA^D$.  Then for any $J \in \tA$,
$[J_{\old}]_{\tA^{\old}} \cong 
[J]_{\tA} \otimes_{\tTT^{\Dp}} (\tTT^{\Dp})^{\old}.$
\item Let $\CA = \CA^{\Dp}$ (resp. $\CA^D$).  Let $\tA$ be
$\tA^{\Dp}$ (resp. $\tA^D$).  Then for any $J \in \CA$,
$[J^4]_{\tA^{\old}} = [J]_{\CA} \otimes_{\TT^{\Dp}} (\tTT^{Dp})^{\old}.$
\end{enumerate}
\end{lemma}
\proof 
We have a natural map
$[J]_{\tA} \rightarrow [J_{\old}]_{\tA^{\old}}$ 
defined by restricting an element
of $[J]$ to its $q_1q_2$-old subvariety.  Locally at any
$\tmm$ outside $S_{\tA}$, this map fits into a commutative
diagram:
$$\begin{CD}
T_{\tmm} \Jmin_{\tA} \otimes [J]_{\tA} @>>> T_{\tmm} J\\
@VVV @VVV\\
T_{\tmm} \Jmin_{\tA^{\old}} \otimes [J^{\old}]_{\tA^{\old}}
@>>> T_{\tmm} J^{\old}.
\end{CD}$$
Since $J^{\old}$ is a subvariety of $J$, the right-hand vertical map
is surjective.  But $T_{\tmm} \Jmin_{\tA^{\old}}$ is free of rank two, so 
the natural map
$([J]_{\tA})_{\tmm} \rightarrow ([J_{\old}]_{\tA^{\old}})_{\tmm}$ 
is surjective as well.  
Statement (1) thus follows by Lemma~\ref{lemma:faithful}. 
Statement (2) is immediate from the previous lemma.
\endproof
 
With these technicalities out of the way, we begin our study of the
geometry of these varieties.  First of all, we have:
\begin{lemma} \label{lemma:twocontrol}
There are non-canonical isomorphisms:
$[\tJ^{D,Dq_1q_2}]_{\tmm} \cong [\tJ^{Dq_1q_2,Dq_1q_2}]_{\tmm}$ and	
$[\tJ^{\Dp,\Dp q_1q_2}]_{\tmm} \cong [\tJ^{\Dp q_1q_2,\Dp q_1q_2}]_{\tmm}.$
\end{lemma}
\proof
We construct these isomorphisms for $\tJ^{D,Dq_1q_2}$ and
$\tJ^{Dq_1q_2,Dq_1q_2}$; the other case is nearly identical.
Since $\tJ^{Dq_1q_2,Dq_1q_2}$ is self-dual, we have
$$[\tJ^{Dq_1q_2,Dq_1q_2}]_{\tmm} \cong [\tJ^{Dq_1q_2,Dq_1q_2}]^*_{\tmm}.$$
Since $\tmm$ is controllable at $q_1$, we have an isomorphism
$$[\tJ^{Dq_1q_2,Dq_1q_2}]_{\tmm} \cong X_{q_1}(\tJ^{Dq_1q_2,Dq_1q_2})_{\tmm}$$
by Proposition~\ref{prop:charisogeny} and Lemma~\ref{lemma:control}.
On the other hand, by Proposition~\ref{prop:chargp}, we have an
isomorphism $$X_{q_1}(\tJ^{Dq_1q_2,Dq_1q_2})^*_{\tmm}
\cong X_{q_2}(\tJ^{D,Dq_1q_2})_{\tmm}.$$  Combining this with the previous
two isomorphisms we find an isomorphism
$$[\tJ^{Dq_1q_2,Dq_1q_2}]_{\tmm} \cong X_{q_2}(\tJ^{D,Dq_1q_2})_{\tmm}.$$
As $\tmm$ is controllable at $q_2$, $X_{q_2}(\tJ^{D,Dq_1q_2})_{\tmm}$ is
isomorphic to $[\tJ^{D,Dq_1q_2}]_{\tmm}$ by Proposition~\ref{prop:charisogeny}
and Lemma~\ref{lemma:control}.
\endproof 

\begin{prop} \label{prop:newfusion} 
We have $\dim \tJ^{Dq_1q_2,Dq_1q_2}[\tmm] = \dim \tJ^{D,D}[\tmm]$.
\end{prop}
\proof
By Lemma~\ref{lemma:twocontrol}, $$\dim \tJ^{Dq_1q_2,Dq_1q_2}[\tmm] =
\dim \tJ^{D,Dq_1q_2}[\tmm].$$  
As $\tJ^{D,Dq_1q_2}$
is an abelian subvariety of $\tJ^{D,D},$
$$\dim \tJ^{D,Dq_1q_2}[\tmm] \leq \dim \tJ^{D,D}[\tmm].$$  
We therefore have
$$\dim \tJ^{Dq_1q_2,Dq_1q_2}[\tmm] \leq \dim \tJ^{D,D}[\tmm],$$ 
so it suffices to show that the opposite inequality holds.

Let $Y_{q_1}$ denote the character group of $\tJ^{Dq_1q_2,Dq_1q_2}[\tmm]$
at $q_1$, and $X_{q_2}$ denote the character group of
$J^D(\Gamma_0(q_2N/D) \cap \Gamma)^2$ at $q_2$.  Since $\tmm$ is
not Eisenstein and has residue characteristic greater than 3,
Proposition~\ref{prop:ribexact2}
implies that $(Y_{q_1}^*/Y_{q_1})_{\tmm}$ is isomorphic to 
$(X_{q_2}/(T_{q_1}^2 - \<q_1\>)X_{q_2})_{\tmm}$.  In particular, we have:
$$\dim Y_{q_1}^*/\tmm Y_{q_1}^* \geq 
\dim Y_{q_1}^*/Y_{q_1} + \tmm Y_{q_1}^* = \dim X_{q_2}/\tmm X_{q_2}.$$

We will relate $Y_{q_1}^*/\tmm Y_{q_1}^*$ to $\tJ^{Dq_1q_2,Dq_1q_2}[\tmm]$,
and $X_{q_2}/\tmm X_{q_2}$ to $\tJ^{D,D}[\tmm]$.  The key is controllability
of $\tmm$ at $q_1$ and $q_2$.

Since $\tmm$ is controllable at $q_1$, $[\tJ^{Dq_1q_2,Dq_1q_2}]_{\tmm}$
is isomorphic to $(Y_{q_1})_{\tmm}$ by Proposition~\ref{prop:charisogeny}
and Lemma~\ref{lemma:control}.
Thus 
$T_{\tmm} \tJ^{Dq_1q_2,Dq_1q_2}$ is
isomorphic to $(Y_{q_1})^2_{\tmm}$.  It follows that
$\tJ^{Dq_1q_2,Dq_1q_2}[\tmm]^{\vee} \cong (Y_{q_1}/\tmm Y_{q_1})^2$.
Since $\tJ^{Dq_1q_2,Dq_1q_2}$ is self-dual,
we also have $[\tJ^{Dq_1q_2,Dq_1q_2}]_{\tmm} \cong
[\tJ^{Dq_1q_2,Dq_1q_2}]^*_{\tmm}$.  In particular $(Y_{q_1})_{\tmm}$ is
also self-dual; that is,
$(Y_{q_1})_{\tmm} \cong (Y^*_{q_1})_{\tmm}$.  We thus find that
$$\dim \tJ^{Dq_1q_2,Dq_1q_2}[\tmm] = 2\dim Y_{q_1}/\tmm Y_{q_1}
= 2\dim (Y_{q_1})^*/\tmm (Y_{q_1})^*.$$

Combining this result with the above inequality, we see that it suffices
to show that $$\dim \tJ^{D,D}[\tmm] = 2\dim X_{q_2}/\tmm X_{q_2}.$$
Since $\tmm$ is controllable at $q_2$,
$(X_{q_2})^2_{\tmm}$ is isomorphic to
$T_{\tmm} J^D(\Gamma_0(q_2N/D) \cap \Gamma)$, by 
Corollary~\ref{cor:controllable}. 
Thus
$$\dim J^D(\Gamma_0(q_2N/D) \cap \Gamma)[\tmm] = 2\dim X_{q_2}/\tmm X_{q_2}.$$
By~\cite{DT}, Theorem 2, the isogeny 
$$J^D(\Gamma_0(q_2N/D) \cap \Gamma)^2 \rightarrow (\tJ^{D,D})_{\qtold}$$ 
has kernel supported away from $\tmm$,
so $\dim (\tJ^{D,D})_{\qtold}[\tmm] = 2\dim X_{q_2}/\tmm X_{q_2}.$
But by Proposition~\ref{prop:fusion}, 
$$\dim (\tJ^{D,D})_{\qtold}[\tmm] = \dim \tJ^{D,D}[\tmm],$$ 
so the result follows.
\endproof

If we fix an isomorphism: $[\tJ^{\Dp,\Dp q_1q_2}]_{\tmm} \rightarrow
[\tJ^{\Dp q_1q_2, \Dp q_1q_2}]_{\tmm}$, we obtain by Theorem~\ref{thm:equiv}
an element of $\Hom(\tJ^{\Dp q_1q_2, \Dp q_1q_2}, \tJ^{\Dp,\Dp q_1q_2})_{\tmm}$.
Composing with the natural inclusion of $\tJ^{\Dp, \Dp q_1q_2}$ into
$\tJ^{\Dp, \Dp}$, we obtain a map $f$ in $\Hom(\tJ^{\Dp q_1q_2, \Dp q_1q_2},
\tJ^{\Dp, \Dp})_{\tmm}$ which induces a surjection $[\tJ^{\Dp, \Dp}]_{\tmm} 
\rightarrow [\tJ^{\Dp q_1q_2, \Dp q_1q_2}]_{\tmm}.$

\begin{lemma}
There is a map $g$ in $\Hom(\tJ^{Dq_1q_2, Dq_1q_2}, \tJ^{D,D})_{\tmm}$
whose dual $g^{\dual}$ fits in a commutative diagram:
$$
\begin{CD}
X_p(\tJ^{Dq_1q_2,Dq_1q_2})_{\tmm} @>>> 
X_{p^{\prime}}(\tJ^{\Dp q_1q_2, \Dp q_1q_2})_{\tmm} \\
@Vg^{\dual}VV @VVf^{\dual}V\\
X_p(\tJ^{D,D})_{\tmm} @>>> 
X_{p^{\prime}}(\tJ^{\Dp,\Dp})_{\tmm},
\end{CD}$$
in which the horizontal maps are those in the exact sequence
of Theorem~\ref{thm:ribexact}, the left-hand vertical map is
induced by $g^{\dual}$, and the right-hand vertical map is
induced by $f^{\dual}$.
\end{lemma}
\proof
As the upper horizontal map in the above diagram identifies
$X_p(\tJ^{Dq_1q_2,Dq_1q_2})_{\tmm}$ with
$X_{p^{\prime}}(\tJ^{D,D})[I]$, where $I$ is the kernel of
the map $\tTT \rightarrow \tTT^{\qqnew}$, and similarly for
the lower horizontal map, it is clear that there is a map
$X_p(\tJ^{Dq_1q_2,Dq_1q_2})_{\tmm} \rightarrow X_p(\tJ^{D,D})_{\tmm}$
making the above diagram commute; the only thing we have to check
is that it is induced by an element $g^{\dual}$ of $\Hom(\tJ^{D,D},
\tJ^{Dq_1q_2, Dq_1q_2})_{\tmm}$.  In particular it suffices to show
that the natural map: $$\Hom((\tJ^{D,D})^{\qqnew},\tJ^{Dq_1q_2, Dq_1q_2})_{\tmm}
\rightarrow \Hom(X_p(\tJ^{Dq_1q_2, Dq_1q_2}), 
X_p((\tJ^{D,D})^{\qqnew}))_{\tmm}$$ is an isomorphism.  (Or, equivalently,
that the natural map $$\Hom([\tJ^{Dq_1q_2, Dq_1q_2}], 
[(\tJ^{D,D})^{\qqnew}])_{\tmm} \rightarrow 
\Hom(X_p(\tJ^{Dq_1q_2,Dq_1q_2}), X_p((\tJ^{D,D})^{\qqnew}))_{\tmm}$$
is an isomorphism.)

By Lemma~\ref{lemma:twocontrol}, the modules
$[\tJ^{D,Dq_1q_2}]_{\tmm}$ and $[\tJ^{Dq_1q_2,Dq_1q_2}]_{\tmm}$
are isomorphic.  It follows that there exists an injection
of $[\tJ^{D,Dq_1q_2}]$ into $[\tJ^{Dq_1q_2,Dq_1q_2}]$ whose cokernel
is supported away from $\tmm$.  By Theorem~\ref{thm:equiv},
this implies that there is an isogeny 
$$\phi:\tJ^{Dq_1q_2,Dq_1q_2} \rightarrow \tJ^{D,Dq_1q_2}$$ 
whose kernel is supported away from $\tmm$.  Then
$\phi^{\dual}$ is an isogeny from
$(\tJ^{Dq_1q_2,Dq_1q_2})^{\dual}$ to $(\tJ^{D,Dq_1q_2})^{\dual}$ whose
kernel is supported away from $\tmm$.  But
$\tJ^{Dq_1q_2, Dq_1q_2}$ is self-dual, and
$$(\tJ^{D,Dq_1q_2})^{\dual} = (\tJ^{D,D}_{\qqnew})^{\dual}
\cong (\tJ^{D,D})^{\qqnew},$$ so we can consider $\phi^{\dual}$
as an isogeny from $(\tJ^{D,D})^{\qqnew}$ to $\tJ^{Dq_1q_2,Dq_1q_2},$
whose kernel is supported away from $\tmm$.
Viewed in this way, $\phi^{\dual}$ induces isomorphisms
between $[\tJ^{Dq_1q_2,Dq_1q_2}]_{\tmm}$ and
$[(\tJ^{D,D})^{\qqnew}]_{\tmm}$, as well as between
$X_p((\tJ^{D,D})^{\qqnew})_{\tmm}$ and
$X_p(\tJ^{Dq_1q_2,Dq_1q_2})_{\tmm}$.

Under these identifications, the natural map:
$$\Hom([\tJ^{Dq_1q_2, Dq_1q_2}],
[(\tJ^{D,D})^{\qqnew}])_{\tmm} \rightarrow
\Hom(X_p(\tJ^{Dq_1q_2,Dq_1q_2}), X_p((\tJ^{D,D})^{\qqnew}))_{\tmm}$$
is identified with the map:
$$\End([\tJ^{Dq_1q_2,Dq_1q_2}]_{\tmm}) \rightarrow
\End(X_p(\tJ^{Dq_1q_2,Dq_1q_2})_{\tmm}).$$
It thus suffices to show that the latter is an isomorphism.
We have a sequence of isomorphisms:
$$\begin{array}{rcl}
\End(X_p(\tJ^{Dq_1q_2,Dq_1q_2})_{\tmm}) & \cong &
\End(X_p(\Jmin_{Dq_1q_2})_{\tmm} \otimes [\tJ^{Dq_1q_2,Dq_1q_2}]_{\tmm}) \\
& \cong & \Hom([\tJ^{Dq_1q_2,Dq_1q_2}]_{\tmm}, \Hom(X_p(\Jmin_{Dq_1q_2}),
X_p(\tJ^{Dq_1q_2,Dq_1q_2}))_{\tmm}) \\
& \cong & \Hom([\tJ^{Dq_1q_2,Dq_1q_2}]_{\tmm}, [\tJ^{Dq_1q_2,Dq_1q_2}]_{\tmm}).
\end{array}$$
(The first isomorphism comes from Proposition~\ref{prop:charisogeny};
the second from the adjointness of $\Hom$ and tensor product, and the
third from Corollary~\ref{cor:key}, using the fact that $\tmm$ is
controllable at $q_1$.)
It follows by Lemma~\ref{lemma:hom} that the natural map in question
is an isomorphism.
\endproof

We are now in a position to prove:
\begin{prop} \label{prop:full}
The natural map $[\tJ^{D,D}]_{\tmm} \rightarrow
\Hom(X_p(\tJmin_D), X_p(\tJ^{D,D}))_{\tmm}$ is an isomorphism.
\end{prop}
\proof
For conciseness of notation, we abbreviate $q_1q_2$ by
$\bq$, and write $X_p$ for
$X_p(\tJmin_D)$ and $\overline{X_p}$ for $X_p(\tJmin_{D\bq})$.
There is a natural surjection $\pi_p: X_p \rightarrow \overline{X_p}$
induced by the inclusion of $\tJmin_{D\bq}$ into $\tJmin_D$.

We have a sequence of isomorphisms:
$$\Hom(X_p, X_p(\tJ^{D,D}))^*_{\tmm} \rightarrow
(X_p)_{\tmm} \otimes X_p(\tJ^{D,D})_{\tmm}^* \rightarrow
(X_p)_{\tmm} \otimes X_{p^{\prime}}(\tJ^{\Dp,D})_{\tmm},$$
where the first isomorphism comes from Lemma~\ref{lemma:star}
and the second from Proposition~\ref{prop:chargp}.  We also obtain
isomorphisms:
$$\Hom(\overline{X_p}, X_p(\tJ^{D\bq,D\bq}))^*_{\tmm} \rightarrow
(\overline{X_p})_{\tmm} \otimes X_p(\tJ^{D\bq,D\bq})_{\tmm}^* 
\rightarrow (\overline{X_p})_{\tmm} \otimes 
X_{p^{\prime}}(\tJ^{\Dp \bq,D\bq})_{\tmm},$$
in the same fashion.  By Lemma~\ref{lemma:compatible} these fit into a
commutative diagram:
$$\begin{CD}
\Hom(X_p, X_p(\tJ^{D,D}))^*_{\tmm} @>>> 
(X_p)_{\tmm} \otimes X_p(\tJ^{D,D})_{\tmm}^* @>>> 
(X_p)_{\tmm} \otimes X_{p^{\prime}}(\tJ^{\Dp,D})_{\tmm} \\
@V(g^{\dual})^*VV @V\pi_p \otimes (g^{\dual})^*VV 
@VV\pi_p \otimes \overline{f}V \\
\Hom(\overline{X_p}, X_p(\tJ^{D\bq,D\bq}))^*_{\tmm} @>>> 
(\overline{X_p})_{\tmm} \otimes X_p(\tJ^{D\bq,D\bq})_{\tmm}^* 
@>>> (\overline{X_p})_{\tmm} \otimes 
X_{p^{\prime}}(\tJ^{\Dp \bq,D\bq})_{\tmm}.
\end{CD}$$

In particular, the right-hand vertical map is surjective, since the
kernel of $f$ is supported away from $\tmm$, so the left-hand
vertical map is surjective as well.  This map fits into a commutative diagram:
$$\begin{CD}
\Hom(X_p, X_p(\tJ^{D,D}))_{\tmm}^* @>>> 
[\tJ^{D,D}]_{\tmm}^*\\
@V(g^{\dual})^*VV @V(g^{\dual})^*VV\\
\Hom(\overline{X_p},X_p(\tJ^{D\bq,D\bq}))_{\tmm}
@>>> [\tJ^{D\bq,D\bq}]_{\tmm}^*,
\end{CD}$$
where the horizontal maps are the $\ZZ_l$-duals of the canonical maps
$$\begin{array}{ccc}
[\tJ^{D,D}] \rightarrow \Hom(X_p,X_p(\tJ^{D,D})) & \text{and} &
[\tJ^{D\bq,D\bq}] \rightarrow \Hom(\overline{X_p}, 
X_p(\tJ^{D\bq,D\bq})).\end{array}$$  
By Corollary~\ref{cor:key}, the bottom horizontal map is an isomorphism. 

Since the left-hand vertical map is a surjection, the right-hand vertical
map is as well.  Moreover, since both $\tJ^{D,D}$ and
$\tJ^{D\bq,D\bq}$ are self-dual, we have $[\tJ^{D,D}]_{\tmm} 
\cong [\tJ^{D,D}]_{\tmm}^*$ and similarly for $[\tJ^{D\bq,D\bq}]$.
By Proposition~\ref{prop:newfusion}, and the fact that
$2 \dim [\tJ^{D,D}] \otimes \tTT/\tmm = \dim \tJ^{D,D}[\tmm]$ (and similarly
for $\tJ^{D\bq,D\bq}$,)
it follows that the right-hand vertical map becomes an isomorphism after 
tensoring with $\tTT/\tmm$.  Thus, after tensoring with $\tTT/\tmm$,
the bottom horizontal map and the right-hand vertical maps are isomorphisms,
and the left-hand vertical map is surjective.  It follows that the
upper horizontal map becomes surjective after tensoring with $\tTT/\tmm$.
By Nakayama's lemma, this implies that the upper horizontal map is
surjective; since both of the modules in the top row have the same
$\ZZ_l$-rank and no $\ZZ_l$-torsion, the upper horizontal map is
an isomorphism.  The result follows immediately.  
\endproof

It remains to pass from $\tJ^{D,D}$ to $(\tJ^{D,D})^{\old}$,
and then to $J^{D,D}$.

\begin{prop} \label{prop:old} 
The natural map:
$[(\tJ^{D,D})^{\old}]_{\tmm} \rightarrow \Hom(X_p(\tJmin_{D,\old})_{\tmm},
X_p((\tJ^{D,D})^{\old})_{\tmm})$ is an isomorphism.
\end{prop}
\proof
Let $\phi$ denote the natural map $\tJ^{D,D} \rightarrow (\tJ^{D,D})^{\old}.$
Composition with $\phi$ induces a natural map
$$[(\tJ^{D,D})^{\old}] = \Hom((\tJ^{D,D})^{\old}, \tJmin_{D,\old}) \rightarrow
 \Hom(\tJ^{D,D}, \tJmin_{D,\old}).$$  Since every map:
$\tJ^{D,D} \rightarrow \tJmin_{D,\old}$ factors through the
$q_1q_2$-old quotient, this map is an isomorphism.  The inclusion
of $\tJmin_{D,\old}$ into $\tJmin_D$ thus induces a map
$[(\tJ^{D,D})^{\old}] \rightarrow [\tJ^{D,D}]$.  Since $\tJmin_{D,\old}$
is the connected component of the identity in the subvariety of
$\tJmin_D$ annilated by the kernel $I$ of the map $\tTT^D \rightarrow
(\tTT^D)^{\old}$, this map identifies $[(\tJ^{D,D})^{\old}]$ with
the submodule of $[\tJ^{D,D}]$ killed by $I$.  It follows that the
dual map $[\tJ^{D,D}] \rightarrow [(\tJ^{D,D})^{\old}]$ is surjective.

We have a commutative diagram:
$$\begin{CD}
\Hom(X_p(\tJmin_D)_{\tmm}, X_p(\tJ^{D,D})_{\tmm})^* @>>> 
[\tJ^{D,D}]_{\tmm}\\
@VVV @VVV \\
\Hom(X_p(\tJmin_{D,\old})_{\tmm}, X_p((\tJ^{D,D})^{\old})_{\tmm})^* 
@>>> [(\tJ^{D,D})^{\old}]_{\tmm},\\
\end{CD}$$
in which the vertical maps are induced by the surjection 
$\tJ^{D,D} \rightarrow (\tJ^{D,D})^{\old}$, and the horizontal maps
are the $\ZZ_l$-duals of the maps taking a morphism of
abelian varieties to the induced maps on character groups.

By the previous proposition, the upper horizontal map is an isomorphism;
the right-hand vertical map is surjective.  It follows that the lower
horizontal map is surjective, and hence an isomorphism.
\endproof

The following technical lemma is now all that we need in order to
establish Theorem~\ref{thm:main} at $\mm$.

\begin{lemma} \label{lemma:four} 
There is an isogeny $(\tJ^{D,D})^{\old} \rightarrow (J^{D,D})^4$
whose kernel is supported away from $\tmm$. 
\end{lemma}
\proof
It is equivalent to show that $[(J^{D,D})^4]_{\tmm} \cong 
[(\tJ^{D,D})^{\old}]_{\tmm}$.  The isogeny $(J^{D,D})^4 \rightarrow
\tJ^{D,D}_{\old}$ has kernel supported away from $\tmm$ by
~\cite{DT}, Theorem 2, so the dual isogeny
$(\tJ^{D,D})^{\old} \rightarrow ((J^{D,D})^4)^{\dual}$ has kernel
supported away from $\tmm$ as well.  In particular,
$$[(\tJ^{D,D})^{\old}]_{\tmm} \cong [((J^{D,D})^4)^{\dual}]_{\tmm}
\cong [(J^{D,D})^4]_{\tmm}^*,$$
so it suffices to show that $[(J^{D,D})^4]_{\tmm}$ is self-dual.

We begin by showing that $(\tTT^D_{\tmm})^* \cong (\TT^D_{\mm})^* 
\otimes \tTT^D_{\tmm}$.  Observe that $\tTT^D$ is free of rank four
over $\TT^D$, with basis $\{1, T_{q_1}, T_{q_2}, T_{q_1}T_{q_2}\}$.
This choice of basis allows us to identify both $(\tTT^D)^*$ and
$(\TT^D)^* \otimes \tTT^D$ with $((\TT^D)^*)^4$ as $\TT^D$-modules.
The Hecke operators $T_{q_1}$ and $T_{q_2}$ act via the matrices
$$\begin{array}{ccc}T_{q_1} = \left( 
\begin{array}{cccc}
0 & - q_1\<q_1\> & 0 & 0\\
1 & \tau_{q_1} & 0 & 0\\
0 & 0 & 0 & -q_1\<q_1\>\\
0 & 0 & 1 & \tau_{q_1}
\end{array} \right); & \quad &
T_{q_2} = \left(
\begin{array}{cccc}
0 & 0 & -q_2\<q_2\> & 0\\
0 & 0 & 0 & -q_2\<q_2\>\\
1 & 0 & \tau_{q_2} & 0\\
0 & 1 & 0 & \tau_{q_2}
\end{array} \right)\end{array}$$
on the latter, and via the transposes of these matrices on the former.
After tensoring with $\ZZ_l$, conjugation with the diagonal matrix
$$P = \left(
\begin{array}{cccc}
q_1\<q_1\>q_2\<q_2\> & 0 & 0 & 0\\
0 & -q_2\<q_2\> & 0 & 0\\
0 & 0 & -q_1\<q_1\> & 0\\
0 & 0 & 0 & 1
\end{array} \right)$$
sends the above matrices to their transposes.  Thus $P$ induces an isomorphism
of $(\TT^D)^* \otimes \tTT^D \otimes \ZZ_l$ with
$(\tTT^D)^* \otimes \ZZ_l$.  Localizing at $\tmm$ then gives the desired
isomorphism.

Now we have:
$$[(J^{D,D})^4]_{\tmm}^* \cong 
\Hom([J^{D,D}]_{\mm} \otimes \tTT^D_{\tmm}, (\tTT^D)^*_{\tmm}) \cong
\Hom([J^{D,D}]_{\mm} \otimes \tTT^D_{\tmm}, (\TT^D)^*_{\mm} \otimes \tTT^D_{\tmm}).$$
The latter is just $[J^{D,D}]_{\mm}^* \otimes \tTT^D_{\tmm}$, which is
isomorphic to $[J^{D,D}]_{\mm} \otimes \tTT^D_{\tmm} 
\cong [(J^{D,D})^4]_{\tmm}$, as required.  \endproof

\noindent{\it Proof of Theorem~\ref{thm:main}}
It suffices to prove that
the natural map
$$[J^{D,D}]_{\mm} \rightarrow \Hom(X_p(\Jmin_D)_{\mm}, X_p(J^{D,D})_{\mm})$$
is an isomorphism.
By Lemma~\ref{lemma:four} and Proposition~\ref{prop:old}, the natural map:
$$[J^{D,D}]^4_{\tmm} \rightarrow \Hom(X_p((\Jmin_D)^4)_{\tmm}, 
X_p((J^{D,D})^4)_{\tmm})$$
is an isomorphism.  It is easy to check that this latter map is
just the corresponding map 
$$[J^{D,D}]_{\mm} \rightarrow \Hom(X_p(\Jmin_D)_{\mm}, X_p(J^{D,D})_{\mm})$$
tensored over $\TT^D_{\mm}$ with $\tTT^D_{\tmm}$.  Since $\tTT^D_{\tmm}$
is free over $\TT^D_{\mm}$, the result is immediate.
\endproof

\section{Main results} \label{sec:main}

Having successfully proven Theorem~\ref{thm:main}, we will now make
use of it to compute $[J^D]$ for $D$ dividing $N$.  Recall that 
for each $D$ dividing $N$ with an even number of prime divisors, we have
constructed a locally-free $\TT^D_{S^D}$-module
$L_D$ such that $[(\Jmin_D)^{\dual}] \cong_{S^D} (\TT^D)^* \otimes L_D$, 
which is unique up to isomorphism (Proposition~\ref{prop:L}).  
We have shown that $[J^{D,D}] \cong_{S^D} [J^{D,D}]^* \otimes L_D$. 

\begin{prop} \label{prop:pq} 
Let $p$ and $q$ be two primes dividing $D$.  Then we have
$$[J^{D,D}] \cong_{S^D} [J^{\frac{D}{pq},D}] \otimes X_p(\Jmin_D) \otimes 
X_q(\Jmin_D) \otimes L_D.$$
\end{prop}
\proof
We have a sequence of isomorphisms:
$$\begin{array}{rcll}
[J^{D,D}] & \cong_{S^D} & [J^{D,D}]^* \otimes L_D &
\mbox{(self-duality of $J^{D,D}$)} \\
& \cong_{S^D} & \Hom(X_p(\Jmin_D),X_p(J^{D,D}))^* \otimes L_D &
\mbox{(Theorem~\ref{thm:main})} \\
& \cong_{S^D} & X_p(\Jmin_D) \otimes X_p(J^{D,D})^* \otimes L_D &
\mbox{(Lemma~\ref{lemma:star})} \\
& \cong_{S^D} & X_p(\Jmin_D) \otimes X_q(J^{\frac{D}{pq},D}) \otimes L_D &
\mbox{(Proposition~\ref{prop:chargp})} \\
& \cong_{S^D} & X_p(\Jmin_D) \otimes X_q(\Jmin_D) \otimes [J^{\frac{D}{pq},D}]
\otimes L_D & \mbox{(Proposition~\ref{prop:charisogeny}),}
\end{array}$$
which immediately proves the claim.
\endproof

To extend this to a relationship between $[J^{D,R}]$ and $[J^{Dpq,R}]$
for $R$ divisible by $Dpq$, we use the following lemma:

\begin{lemma} \label{lemma:quotient1}
Let $D$ $R$, and $R^{\prime}$ be divisors of $N$ such that
$D$ is divisible by an even number of
primes, and such that $D | R | R^{\prime}$. 
Then $[J^{D,R^{\prime}}] \cong_{S^{R^{\prime}}} [J^{D,R}] 
\otimes_{\TT^R} \TT^{R^{\prime}}$.
\end{lemma}
\proof
There is a natural map $[J^{D,R}] \rightarrow [J^{D,R^{\prime}}]$ given by
taking an element of $[J^{D,R}]$, treating it
as an element of $\Hom(J^{D,R}, \Jmin_R),$ and
restricting it to $J^{D,R^{\prime}} \subset J^{D,R}$.  (The image of 
$J^{D,R^{\prime}}$ will
always land in $\Jmin_{R^{\prime}} \subset \Jmin_R$.)  
By Lemma~\ref{lemma:faithful},
it suffices to show the surjectivity of this map at $\mm$ outside
$S_{R^{\prime}}$.  Let $\hat \mm$ be the preimage of $\mm$ in $\TT^R$.
We have a commutative diagram:
$$\begin{CD}
T_{\hat \mm} \Jmin_R \otimes [J^{D,R}] @>>> 
T_{\mm} \Jmin_{R^{\prime}} \otimes [J^{D,R^{\prime}}]_{\mm}\\
@VVV @VVV \\
T_{\hat \mm} J^{D,R} @>>> T_{\mm} J^{D,R^{\prime}},
\end{CD}$$
in which the vertical maps are isomorphisms.  The lower horizontal map
is induced by an inclusion of varieties and is therefore surjective.
The upper horizontal map is the restriction map we wish to study, tensored
with $T_{\hat \mm} \Jmin_R$.  The surjectivity of the lower horizontal map
means that the upper horizontal map is surjective as well.  But since
$T_{\hat \mm} \Jmin_R$ is free of rank two over $\TT^R_{\mm}$, this implies 
that the map $[J^{D,R}]_{\hat \mm} \rightarrow [J^{D,R^{\prime}}]_{\mm}$ 
is surjective as well.
\endproof

We also need to relate $L_R$ with $L_{R^{\prime}}$ for $R$
dividing $R^{\prime}$; since we have chosen our $\Jmin_R$ compatibly,
we have the following:

\begin{lemma} \label{lemma:quotient2}
If $R$ divides $R^{\prime}$, then $L_{R^{\prime}}$ is 
$S^{R^{\prime}}$-isomorphic to $L_R \otimes_{\TT^R} \TT^{R^{\prime}}$.
\end{lemma}
\proof The property 
$$[(\Jmin_{R^{\prime}})^{\dual}] 
\cong_{S^{R^{\prime}}} (\TT^{R^{\prime}})^*
\otimes L_{R^{\prime}}$$ 
characterises $R^{\prime}$ up to $S^{R^{\prime}}$-isomorphism.
Since $\Jmin_{R^{\prime}}$ is the $\frac{R^{\prime}}{R}$-new
subvariety of $\Jmin_R$, $\Jmin_{R^{\prime}}$ is the
connected component of $\Jmin_R[I],$ where $I$ is the
kernel of the map $\TT^R \rightarrow \TT^{R^{\prime}}$.  Thus
$[(\Jmin_{R^{\prime}})^{\dual}] = \Hom((\Jmin_{R^{\prime}})^{\dual},
\Jmin_{R^{\prime}})$ is isomorphic to
$\Hom((\Jmin_{R^{\prime}})^{\dual}, \Jmin_R)[I]$.  The
injection $\Jmin_{R^{\prime}} \rightarrow \Jmin_R$ dualizes to a surjection
$(\Jmin_R)^{\dual} \rightarrow (\Jmin_{R^{\prime}})^{\dual}$; moreover,
any morphism of $(\Jmin_R)^{\dual}$ that is killed by $I$ factors through
$(\Jmin_{R^{\prime}})^{\dual}$.  Thus 
$$\Hom((\Jmin_{R^{\prime}})^{\dual}, \Jmin_R)[I] \cong_{S^{R^{\prime}}}
 \Hom((\Jmin_R)^{\dual}, \Jmin_R)[I] = [(\Jmin_R)^{\dual}][I].$$

It follows that $[(\Jmin_{R^{\prime}})^{\dual}] \cong [(\Jmin_R)^{\dual}][I]$.
Dualizing, and invoking Corollary~\ref{cor:faithful}, we find
that $L_{R^{\prime}}^{-1} \cong L_R \otimes_{\TT^R} \TT^{R^{\prime}}$,
as required.
\endproof

Putting these together, and passing to $\CA^N$, we find:
\begin{prop}
Let $D$ divide $N$ and be divisible by an even number of primes.  For
any $p$ and $q$ dividing $D$, we have:
$$[J^D] \cong_{S^D} [J^{\frac{D}{pq}}] \otimes X_p(\Jmin) \otimes X_q(\Jmin)
\otimes L_N.$$
\end{prop}
\proof Tensor the isomorphism of Proposition~\ref{prop:pq} with $\TT^N$
and apply Lemmas~\ref{lemma:quotient1} and~\ref{lemma:quotient2}. 
\endproof

Working inductively, one immediately proves:
\begin{thm} \label{thm:global1}
For any $D$ dividing $N$ and divisible by $2k$ primes, we have:
$$[J^D] \cong_S [J^1] \otimes (L_N)^{\otimes k} \otimes \bigotimes_{p|D}
X_p(\Jmin).$$
\end{thm}

Combining this with Proposition~\ref{prop:JJ'}, we immediately
obtain:
\begin{cor} \label{cor:global1}
Let $D$ and $D^{\prime}$ be divisors of $N$, divisible by $2n$ and
$2m$ primes, respectively, such that $D$ divides $D^{\prime}$.
Let $M$ be the module  
$$(L_N)^{\otimes m-n} \otimes \bigotimes_{p|\frac{D^{\prime}}{D}}
X_p(\Jmin).$$
Then there is an $S$-isomorphism: 
$$g:\Hom(J^D,J^{D^{\prime}}) \rightarrow \Hom(M,\End(J^D))$$ 
such that if $\phi \in \Hom_S(J^D,J^{D^{\prime}})$,
and $I_{\phi} \subset \End(J^D)$ is the ideal generated by image of $M$ under
$g(\phi): M \rightarrow \End(J^D),$ then $\ker_S \phi = J^D[I_{\phi}]$
up to support on $S$.
\end{cor} 

This result determines the $S$-isomorphism class of $J^{D^{\prime}}$
in terms of $J^D$; in particular, if one has $J^D$, and one wishes to
construct an abelian variety $S$-isomorphic to $J^{D^{\prime}}$, one
merely has to find an ideal $I$ of $\TT$ (or, more generally,
of $\End(J_D)$) isomorphic to $M$, and construct $J^D/J^D[I]$.
Unfortunately, although $M$ is independent of the choice of $\Jmin$,
the various factors in the product defining $M$ are not, making $M$
difficult to compute.  One can remedy this, however, by
introducing a suitable ``multiplicity one'' hypothesis.

In particular, we expand our set $S$ to include those maximal ideals
for which $J^{1,1} = J^1(\Gamma_0(N) \cap \Gamma)$ fails to satisfy
``multiplicity one''; that is, we take $S^R_+$ to be the set of
$\mm \subset \TT^R$ such that either $\mm$ is in $S^R$ or
$J^1(\Gamma_0(N) \cap \Gamma)[\mm]$ has dimension greater than two.
(The latter is expected to be a very rare occurrence; in particular it
cannot occur if $l$ is prime to $2NM$~\cite{Ti}.  Moreover, the only known 
counterexamples are in residue characteristic $2$~\cite{Kilford}.)
As usual we take $S_+ = S^N_+$.

The advantage of enlarging our set of ``bad primes'' to $S_+$ is that
we can take $\Jmin_R = J^{1,R}$ for all $R$; the ``multiplicity one''
condition we have imposed above makes this a valid choice.  Moreover,
we can compute $L_N$ in this setting.  Since $J^{1,1}_{\chi} \cong
(J^{1,1})^{\vee}$, and $[J^{1,1}]_{\chi} \cong_{S^1_+} \TT^1$, we have
$\TT^1 \cong_{S^1_+} L^1 \otimes (\TT^1)^*$.  We thus recover the
result that $\TT^1_{S^1_+}$ is Gorenstein, and find that $L^1 \cong_{S^1_+}
((\TT^1_{S^1_+})^*)^{-1}.$  Lemma~\ref{lemma:quotient2} then determines
$L_N$ explicitly.  Putting this together, we find:

\begin{thm} \label{thm:globalmult1}
For any $D$ dividing $N$ and divisible by exactly $2k$ primes, we have:
$$[J^D] \cong_{S_+} (L_N)^{\otimes k} \otimes \bigotimes_{p|D}
X_p(J^1).$$
\end{thm}

\begin{cor} \label{cor:globalmult1}
Let $D$ and $D^{\prime}$ be divisors of $N$, divisible by exactly $2n$ and
$2m$ primes, respectively, such that $D$ divides $D^{\prime}$.
Let $M$ be the module  
$$(L_N)^{\otimes m-n} \otimes \bigotimes_{p|\frac{D^{\prime}}{D}}
X_p(J^1).$$
Then there is an $S_+$-isomorphism: $$g:\Hom(J^D,J^{D^{\prime}}) \rightarrow
\Hom(M,\End(J^D))$$ 
such that if $\phi \in \Hom_{S_+}(J^D,J^{D^{\prime}})$,
and $I_{\phi} \subset \End(J^D)$ is the ideal generated by image of $M$ under
$g(\phi): M \rightarrow \End(J^D),$ then $\ker_{S_+} \phi = J^D[I_{\phi}]$
up to support on $S_+$.
\end{cor} 

\begin{remark} \rm
The previous two results, although valid over a smaller
subset of $\spec \TT$,
have the advantage that the character groups that appear within them are 
explicitly computable, via the techniques in~\cite{Kohel} together with
Corollary~\ref{cor:newchargp}.  It is thus possible, given the
isomorphism class of $J^1$, to {\em effectively compute} an abelian variety
which is $S_+$-isomorphic to $J^D$ for any $D$ one desires. 
\end{remark}

The following result, a strengthening of a result proved by Yang~\cite{Yang}
in the case in which $\overline{\rho}_{\mm}$ is ramified at at least half of the
primes dividing $D$, is an immediate corollary.

\begin{cor} \label{cor:multbound}
Let $\mm$ be a maximal ideal of $\TT$ outside $S$, and
let $k$ be the number of primes dividing $D$ at which $\mm$ is not
controllable.  Then $\dim J^D[\mm] \leq 2^k \dim J^1[\mm]$.
\end{cor}
\proof
Since $J^D[\mm]$ is naturally isomorphic to 
$(T_{\mm} J^D/\mm T_{\mm} J^D)^{\vee}$, and
similarly for $J_1[\mm]$, it suffices to show that 
$$\dim T_{\mm} J^D/\mm T_{\mm} J^D \leq 2^k \dim J^1/\mm T_{\mm} J^1.$$
Moreover, we have isomorphisms
$$T_{\mm} J^D \cong [J^D]_{\mm} \otimes T_{\mm} \Jmin,$$
and similarly for $T_{\mm} J^1$, so it suffices to show that
$$\dim [J^D]/\mm [J^D] \leq 2^k \dim [J^1]/\mm [J^1].$$
By Nakayama's lemma, this is equivalent to showing that the size of
a minimal generating set for $[J^D]_{\mm}$ is at most $2^k$ times
the size of a minimal generating set for $[J^1]_{\mm}$.

By Theorem~\ref{thm:global1}, we have 
$$[J^D]_{\mm} \cong [J^1]_{\mm} \otimes \bigotimes_{p|D} X_p(\Jmin)_{\mm}.$$
If $\mm$ is controllable at $p$, $X_p(\Jmin)_{\mm}$
is free of rank one, and can be ignored.  If $\mm$ is not controllable
at $p$, the surjection $T_{\mm} \Jmin \rightarrow X_p(\Jmin)_{\mm}$
shows that $X_p(\Jmin)/ \mm X_p(\Jmin)$ has dimension at most two.  In
particular, $X_p(\Jmin)_{\mm}$ is generated by a two-element set, so
tensoring with $X_p(\Jmin)_{\mm}$ increases the size of a minimal
generating set by at most a factor of two.  The result follows. 
\endproof

Considering the special case in which $\mm$ is controllable at
{\em every} prime dividing $D$ yields a ``multiplicity one''
result for Jacobians of Shimura curves:

\begin{cor} \label{cor:mult1}
Let $\mm$ be a maximal ideal of $\TT$ lying outside $S$.
If $\mm$ is controllable at every prime dividing $D$, then
$J^D[\mm]$ and $J^1[\mm]$ have the same dimension.  In particular, if 
$\mm$ satisfies
``multiplicity one'' for $J^{1,1} = J^1(\Gamma_0(N) \cap \Gamma)$, then
$\dim J^D[\mm] = 2.$
\end{cor}
\proof By Theorem~\ref{thm:global1}, and the fact that
$X_p(\Jmin)_{\mm}$ is free for every $p$ dividing $D$, we have
$[J^D]_{\mm} \cong [J^1]_{\mm}$.  It follows that 
$T_{\mm} J^D \cong T_{\mm} J^1,$ and hence that $J^D[\mm] \cong J^1[\mm]$.
If $\mm$ satisfies ``multiplicity one'' for the full Jacobian
$J^{1,1}$,
then it also satisfies ``multiplicity one'' for 
$J^1$, as the latter is just the $N$-new subvariety of the former.
In particular $J^1[\mm]$ has dimension two and the result follows.
\endproof

\begin{remark} \rm Although we have been working throughout with maximal
ideals $\mm$ outside $S$ (and hence of residue characteristic $2$ or $3$), the
preceding corollary holds in more generality.  In particular, if
$\Gamma = \Gamma_0(M)$ or if $\Gamma$ contains $\Gamma_1(r)$, 
for $r \geq 4$, then the
preceding corollary holds for all non-Eisenstein $\mm$.  To establish
this, one observes observe that the only places in the above argument
in which we need the maximal ideals of residue characteristic 
$2$ or $3$ to lie in
$S$ are in Proposition~\ref{prop:component} and in the level-raising
arguments of Section~\ref{sec:raise}.  Since the latter are unnecessary
for an $\mm$ which is already controllable at every prime dividing $D$,
and Proposition~\ref{prop:component} holds even in residue characteristics $2$
and $3$ when $\Gamma = \Gamma_0(M)$ or $\Gamma_1(M)$~\cite{ribchar},
all of the results we have obtained above hold locally at $\mm$.  This
is enough to establish the ``multiplicity one'' result above.
\end{remark}

We also have the following alternative characterization of the
relationship between $[J^D]$ and $[J^{D^{\prime}}]$ for
$D$ dividing $D^{\prime}$: 

\begin{prop} \label{prop:homDD'}
Let $D$ and $D^{\prime}$ be divisors of $N$, each divisible by
an even number of primes, and suppose that $D$ divides $D^{\prime}$.
Then $$\Hom(J^{D^{\prime},D^{\prime}}, J^{D,D^{\prime}})
\cong_{S^{D^{\prime}}} 
\Hom(X_p(J^{D,D^{\prime}}), X_p(J^{D^{\prime},D^{\prime}}))$$
for any $p$ dividing $D^{\prime}$.
\end{prop}
\proof We have an isomorphism $$X_p(J^{D,D^{\prime}}) 
\cong_{S^{D^{\prime}}}
X_p(\Jmin_{D^{\prime}}) \otimes [J^{D,D^{\prime}}].$$  This induces
an isomorphism
$$\Hom(X_p(J^{D,D^{\prime}}),X_p(J^{D^{\prime},D^{\prime}})) 
\cong_{S^{D^{\prime}}} \Hom([J^{D,D^{\prime}}], \Hom(X_p(\Jmin_{D^{\prime}}),
X_p(J^{D^{\prime},D^{\prime}}))),$$ by 
the adjointness of $\Hom$ and tensor products.
By Theorem~\ref{thm:main}, the latter module is naturally 
$S^{D^{\prime}}$-isomorphic to
$\Hom([J^{D,D^{\prime}}], [J^{D^{\prime},D^{\prime}}]),$ and hence
to $\Hom(J^{D^{\prime},D^{\prime}},J^{D,D^{\prime}})$ by 
Theorem~\ref{thm:equiv}.\endproof

\begin{cor} \label{cor:homDD'}
We have a natural isomorphism: $$[J^{D^{\prime},D^{\prime}}] 
\cong_{S^{D^{\prime}}} [J^{D,D^{\prime}}] \otimes 
\Hom(X_p(J^{D,D^{\prime}}),X_p(J^{D^{\prime},D^{\prime}})).$$
\end{cor}
\proof
By Theorem~\ref{thm:global1}, there is an $M$ such that
$[J^{D^{\prime},D^{\prime}}] \cong_{S^{D^{\prime}}} M \otimes 
[J^{D,D^{\prime}}].$  By Lemma~\ref{lemma:hom}, it follows that
$[J^{D^{\prime},D^{\prime}}] \cong_{S^{D^{\prime}}}
\Hom([J^{D,D^{\prime}}],[J^{D^{\prime},D^{\prime}}]) \otimes
[J^{D,D^{\prime}}]$.  The result follows by Proposition~\ref{prop:homDD'}.
\endproof

\begin{cor} \label{cor:homtwo} 
Fix $D$ and $D^{\prime}$ as above, and let $p$ divide
$D^{\prime}$.  Let $M$ be the module 
$$\Hom(X_p(J^{D,D^{\prime}}),X_p(J^{D^{\prime},D^{\prime}})) 
\otimes_{\TT^{D^{\prime}}} \TT.$$  Then $[J^{D^{\prime}}] = M \otimes [J^D]$,
and there is an $S$-isomorphism: $$g:\Hom(J^D,J^{D^{\prime}}) \rightarrow
\Hom(M,\End(J^D))$$ 
such that if $\phi \in \Hom_S(J^D,J^{D^{\prime}})$,
and $I_{\phi} \subset \End(J^D)$ is the ideal generated by image of $M$ under
$g(\phi): M \rightarrow \End(J^D),$ then $\ker_S \phi = J^D[I_{\phi}]$
up to support on $S$.
\end{cor}
\proof Immediate from Corollary~\ref{cor:homDD'}, Lemma~\ref{lemma:quotient1},
and Proposition~\ref{prop:JJ'}.\endproof

\section{Further questions}

These results leave several questions unanswered.  The most obvious of
these is whether or not the set $S$ can be made smaller.  At the
moment, a maximal ideal $\mm$ can be in $S$ either because it is
Eisenstein or because it has residue characteristic $2$ or $3$.  

It seems likely that in many cases the non-Eisenstein maximal ideals of 
residue characteristic $2$ and $3$ may be removed from $S$.  At 
the moment there
are two places in the argument that require the residue characteristic
to be greater than $3$.  The first is that component groups of
Jacobians of some Shimura curves can have support at non-Eisenstein primes
of residue characteristic $2$ or $3$, which interferes with the proof of
Proposition~\ref{prop:chargp}.  This does not happen, however, if
$\Gamma$ either has the form $\Gamma_0(M)$ for some $M$, or if $\Gamma$
contains $\Gamma_1(r)$ for some $r \geq 4$~\cite{ribchar}.  Thus in
many cases this does not pose a problem.

The more serious difficulty occurs in Section~\ref{sec:raise},
where we rely on results of~\cite{DT} which only hold in residue characteristic
greater than $3$.  The only result which we really need from that
paper is Theorem 2 (we invoke Theorem A at one point but if one has Theorem 
2 then the case we need of Theorem A seems to follow quickly.)  Moreover,
since we have the luxury of choosing the primes at which we raise the level,
we need only find two primes $q_1$ and $q_2$ for which Theorem 2 holds
locally at $\mm$.  It seems very unlikely that this will not happen.  It
is thus reasonable to conjecture that (when $\Gamma$ satisfies one of
the conditions given above) Theorem~\ref{thm:main} holds even when $S$ contains
only Eisenstein primes. 
 
The non-Eisenstein condition, on the other hand, is much more serious.  
All of the techniques of section~\ref{sec:isogeny}, and 
Proposition~\ref{prop:chargp}, fail to hold at Eisenstein maximal ideals.
Thus applying the techniques of this paper in this case seems hopeless.

We have also left open the question of whether or not there exist
{\em canonical} maps between the varieties in question.  If one is
interested in looking for such things, one place to start would be to
try to construct canonical maps between the modules constructed in
Section~\ref{sec:main}.

Finally, Corollary~\ref{cor:mult1} provides a sufficient condition
for ``multiplicity one'' for Jacobians of Shimura curves.  It is
interesting to ask if this condition is necessary as well.  This
would hold if one could establish a converse to Lemma~\ref{lemma:control}.
A special case of such a result appears in~\cite{Israel}, where
Ribet constructs character groups which are not free of rank one at
specific maximal ideals of the Hecke algebra. 

\textsc{Acknowledgements.} The author would like to thank Kenneth A.
Ribet for bringing the problem to his attention, and for his continued
advice and encouragement.  The author was partially supported by the 
National Science Foundation. 

\providecommand{\bysame}{\leavevmode\hbox to3em{\hrulefill}\thinspace}

\section{Appendix: The Deligne-Rapoport Model for Shimura curves}

Several of the proofs in section~\ref{sec:chargp} rely on the existence
of a certain model for Shimura curves over a prime of bad reduction.  
Buzzard~\cite{Buz:DR} constructed such a model in the case where the
Shimura curve in question is a fine moduli space; i.e., when the congruence
subgroup $\Gamma$ giving the level structure is sufficiently small.  Here
we extend Buzzard's results to include the case when the Shimura curve is
only a coarse moduli space. 

Let $D$ be a squarefree product of an even number of primes, and let $B$
be an indefinite quaternion algebra over $\QQ$ of discriminant $D$.  Let
$\OO$ be a maximal order of $B$.  Fix
a congruence subgroup $\Gamma$ of level prime to $D$, and consider the
Shimura curves $X = X^D(\Gamma),$ and $X_p = X^D(\Gamma_0(p) \cap \Gamma)$.

We first assume that $\Gamma$ is contained in $\Gamma_1(r)$ for some $r > 3$.
In this case (see for instance~\cite{Buz:DR}) the curves $X$ and
$X_p$ are fine moduli spaces for ``false elliptic curves'' for $B$ with level
structure.  More precisely:

\begin{defn} A false elliptic curve for $\OO$ is an abelian surface $A$,
together with an action of $\OO$ on $A$, such that for all 
$\tau \in \OO$, the reduced trace of $\tau$ is equal to the trace
of $\tau$ acting on the tangent space to $A$ at the identity.
\end{defn}

Then $X$ is the fine moduli space paramaterizing false elliptic curves
for $\OO$ with $\Gamma$-level structure, and $X_p$ is the fine moduli space 
parametrizing triples $(A,P,\rho)$, where $A$ is a false elliptic curve
for $\OO$, $P$ is a $\Gamma$-level structure on $A$, and 
$\rho: A \rightarrow B$ is an $\OO$-equivariant isogeny of false elliptic
curves, of degree $p^2$.

When we restrict to the fiber over $\overline{\FF}_p$, we obtain two
morphisms 
$\Frob, \Ver: X_{\overline{\FF}_p} \rightarrow (X_p)_{\overline{\FF}_p}.$
The former takes a pair $(A,P)$ to the triple 
$(A, P, \Frob: A \rightarrow A^{(p)})$, and the latter takes $(A,P)$
to the triple $(A^{(p)}, P^{(p)}, \Ver: A^{(p)} \rightarrow A)$.  
Then one has:
\begin{thm} (Buzzard, \cite{Buz:DR}) 
\begin{enumerate}
\item $X$ and $X_p$ are regular schemes, defined over $\ZZ[\frac{1}{ND}]$.
\item $X$ is smooth over $\ZZ[\frac{1}{ND}]$, and $X_p$ is smooth over
$\ZZ[\frac{1}{NDp}]$.
\item The maps $\Frob$ and $\Ver$ are closed immersions, and their images
are the two irreducible components of $(X_p)_{\overline{\FF}_p}$.  These
two components intersect transversely at the supersingular points of 
$(X_p)_{\overline{\FF}_p}$.
\item At a supersingular point $x$ of $(X_p)_{\overline{\FF}_p}$, the completion
of the strict henselization of the local ring is isomorphic to
$W(\overline{\FF}_p)[[X,Y]]/(XY-p).$
\end{enumerate}
\end{thm}

In the case where $\Gamma$ is arbitrary, $X$ and $X_p$ need not be
fine moduli spaces.  Let $\CX$ (resp. $\CX_p$) denote the moduli stack
for pairs $(A,P)$ where $A$ is a false elliptic curve for $\OO$ and
$P$ is a $\Gamma$-level structure on $A$ (resp. triples $(A,P,\rho)$ where
$A$ and $P$ are as before and $\rho$ is an $\OO$-equivariant isogeny
of degree $p^2$).  Then we can define maps 
$\Frob, \Ver: \CX_{\overline{\FF}_p} \rightarrow (\CX_p)_{\overline{\FF}_p}$
exactly as before.  These maps induce maps $\Frob, \Ver: X \rightarrow X_p$,
as $X$ and $X_p$ are the underlying coarse moduli spaces of $\CX$ and $\CX_p$.
We then have:

\begin{thm} For arbitrary congruence subgroups $\Gamma$, one has:
\begin{enumerate} 
\item $X$ is a regular scheme over $\ZZ[\frac{1}{ND}]$; $X_p$ is defined
over $\ZZ[\frac{1}{ND}]$ and regular away from the supersingular points on
the fiber over $\overline{\FF}_p$. 
\item $X$ is smooth over $\ZZ[\frac{1}{ND}]$, and $X_p$ is smooth over
$\ZZ[\frac{1}{NDp}]$.
\item The maps $\Frob$ and $\Ver$ are closed immersions, and their images
are the two irreducible components of $(X_p)_{\overline{\FF}_p}$.  These
two components intersect transversely at the supersingular points of
$(X_p)_{\overline{\FF}_p}$,
\item At a supersingular point $x$ of $(X_p)_{\overline{\FF}_p}$, the
completion of the strict henselation of the local ring is isomorphic to
$W(\overline{\FF}_p)[[X,Y]]/(XY-p^{k_x})$, where $k_x$ is the order
of the automorphism group of the false elliptic curve with level structure
corresponding the the point $x$ modulo $\{\pm1\}$ if $-1$ is in $\Gamma$,
and the order of the full automorphism group if $-1$ is not in $\Gamma$.
\end{enumerate}
\end{thm}
\proof
Statement 2 is well-known, and statement 1 is known except for the
regularity of $X_p$ at non-supersingular points on the characteristic $p$
fiber.  It thus suffices to prove this regularity, along with statements
3 and 4.  The arguments we give are basically adaptations of arguments found 
in~\cite{DR} for the case when $D=1$.
Let $r > 3$ be a prime which does not divide $NDp$, and consider
the curves $X_r = X^D(\Gamma \cap \Gamma_1(r))$ and $X_{p,r} =
X^D(\Gamma_0(p) \cap \Gamma \cap \Gamma_1(r))$.  Then $X_r$ and $X_{p,r}$
are fine moduli spaces, so the above theorem applies.

We have maps $X_r \rightarrow X$ and $X_{p,r} \rightarrow X_p$ obtained
by forgetting the level structure at $r$.  These maps are finite, and
ramified precisely at those $x$ in $X$ or $X_p$ for which $k_x > 1$.  
Moreover, 
at any $x$ in $X$ (resp. $X_p$) for which $k_x > 1$, the automorphism
group of $x$ acts on the local ring of any preimage of $x$ in $X_r$
(resp. $X_{p,r}$).  This action is faithful if $-1$ is not in $\Gamma$,
and has kernel $\{\pm 1\}$ if $-1$ is in $\Gamma$.  Let $G_x$ denote the
group $\Aut(x)$ in the former case, and $\Aut(x)/\{\pm 1\}$ in the latter.
Then the above maps induce isomorphisms:
$$\hat \OO_{X,x} \cong \hat \OO_{X_r,\tx}^{G_x}$$
$$\hat \OO_{X_p,x} \cong \hat \OO_{X_{p,r},\tx}^{G_x},$$
where $\tx$ is any lift of $x$ to $X_r$ or $X_{p,r}$, respectively,
and $\hat \OO_{X,x}$ is the completion of the strict henselization of the
local ring of $X$ at $x$.

Suppose $x$ is a point on the characteristic $p$ fiber of $X_p$.
If $\tx$ is a regular point of $X_{p,r}$, then
$\hat \OO_{X_{p,r},\tx}$ is isomorphic to $W(\overline{\FF}_p)[[X]]$.
By the same argument as the proof of~\cite{DR}, VI.6.9, the subring
of $G_x$-invariants of this ring is the ring $W(\overline{\FF}_p)[[X^{k_x}]].$
In particular, $x$ is a regular point.  This proves statement 1.

If $\tx$ is a supersingular point of $X_{p,r}$,
then $\hat \OO_{X_{p,r},\tx}$ is isomorphic to 
$W(\overline{\FF}_p)[[X,Y]]/(XY-p)$.  Again, the same argument as the proof
of~\cite{DR}, VI.6.9 shows that the subring of $G_x$-invariants is
generated by $X^{k_x}$ and $Y^{k_x}$, and hence isomorphic to
$W(\overline{\FF}_p)[[X,Y]]/(XY-p^{k_x})$.  This proves statement 4.

For statement 3, we have a commutative diagram
$$\begin{CD}
(X_r)_{\overline{\FF}_p} @>>> (X_{p,r})_{\overline{\FF}_p} \\ 
@VVV @VVV \\
X_{\overline{\FF}_p} @>>> (X_p)_{\overline{\FF}_p} 
\end{CD}$$
in which the horizontal maps can be either both $\Frob$ or both $\Ver$.
Since the map from $X_{p,r}$ to $X_p$ is a surjection, $X_p$ will thus
be the union of the images of $\Frob$ and $\Ver$, and these images will
intersect precisely at the supersingular points of $X_p$.

All that remains is to prove that $\Frob$ and $\Ver$ are closed immersions.
They are clearly injections on closed points.  Thus, since $X$ is
proper, to prove they are closed immersions it suffices to show that they 
are injections on tangent spaces.
First observe that for any $x \in X$, $G_x = G_{\Frob(x)}$, since any 
automorphism of a false elliptic curve commutes with Frobenius. 
Thus if we fix a lift $\tx$ of $x$ to $X_r$, then $G_x$ acts on 
both $\hat \OO_{X_r, \tx}$ and $\hat \OO_{X_{p,r}, \Frob(\tx)}$.
Moreover, the map from $\hat \OO_{X_p, \Frob(x)}$ to $\hat \OO_{X,x}$
induced by $\Frob$ is just the restriction of the corresponding map 
from $\hat \OO_{X_{p,r}, \Frob(\tx)}$ to $\hat \OO_{X_r, \tx}$ to the 
ring of $G_x$-invariants of $\hat \OO_{X_{p,r}, \Frob(\tx)}$, by the 
commutative diagram given above.  In particular, our explicit computation
of these rings shows that it is surjective.  Thus $\Frob$ induces an
injection on tangent spaces, and is therefore a closed immersion.  The same
argument shows that $\Ver$ is a closed immersion.\endproof

\end{document}